\documentclass{amsart}
\usepackage{graphicx}
\usepackage{amsmath}
\usepackage{amsfonts}
\usepackage{amssymb}
\usepackage{ifpdf}
\newtheorem{theorem}{Theorem}

\newtheorem{corollary}[theorem]{Corollary}

\newtheorem{lemma}[theorem]{Lemma}

\newtheorem{proposition}[theorem]{Proposition}
\newtheorem{remark}[theorem]{Remark}

\newcommand{\dint}{\displaystyle\int}

\newcommand{\dsup}{\displaystyle\sup}

\begin{document}

\title[]{Non-collapsing condition and Sobolev embeddings for  Haj{\l}asz-Besov spaces}
\author{Joaquim Mart\'{i}n$^{\ast}$}
\address{Department of Mathematics\\
Universitat Aut\`{o}noma de Barcelona} \email{jmartin@mat.uab.cat}
\author{Walter A.  Ortiz**}
\address{Department of Mathematics\\
Universitat Aut\`{o}noma de Barcelona} \email{waortiz@mat.uab.cat}
\thanks{$^{\ast}$Partially supported by Grants PID2020-113048GB-I00 and PID2020-114167GB-I00 funded  both by MCIN/AEI/10.13039/501100011033}
\thanks{**Partially supported by Grant 2017SGR395 (AGAUR, Generalitat de Catalunya)}
\thanks{This paper is in final form and no version of it will be submitted for
publication elsewhere.} 
\thanks{Conflict of Interest: The authors declare that they have no conflict
of interest.}

\begin{abstract}
In this paper we will focus on understanding the relation between
Sobolev embedding theorems for Haj{\l}asz-Besov spaces defined on a
doubling metric measure space $(\Omega,d,\mu)$ and the
non-collapsing condition of the measure, i.e.
\[
\inf_{x\in\Omega}\mu(B(x,1))>0.
\]
We will also obtain embedding results for Haj{\l}asz-Besov spaces
whose modulus of smoothness is generated by a rearrangement
invariant quasi-norm.
\end{abstract}
\maketitle

\section{Introduction}

In recent years, an important area of researches in metric measure
spaces has been focused on understanding the relation between
Sobolev embedding theorems and growth properties of the underlying
measure (see for example \cite{Ka01}, \cite{AYY}, \cite{AGH},
\cite{Gor}, \cite{Ha2}, \cite{HKa}, \cite{MaOr}, \cite{MaOr1} and
the references quoted therein). In this work we will continue with
this topic for Haj{\l}asz-Besov spaces.

Throughout the paper $\left( \Omega ,d,\mu \right) $ denotes a
doubling metric measure space, i.e. we assume that the measure $\mu
$ satisfies the following doubling condition: there exists a
constant $C_{\mu }>1$, such that
\begin{equation}
0<\mu (B(x,2r))\leq C\mu (B(x,r))<\infty  \label{dobla1}
\end{equation}
for every ball $B(x,r)$, for all $x\in \Omega $ and $r>0.$ The
smallest constant $C$ for which (\ref{dobla1}) is satisfied will be
called the doubling constant and will be denoted by $C_{\mu }$. The
above doubling condition is equivalent to that, for any ball $B$ and
any $\lambda \in
\lbrack 1,\infty ),$%
\begin{equation*}
\mu (B(x,\lambda r))\leq C_{\mu }\lambda ^{Q}\mu (B(x,r))
\end{equation*}
where $Q:=\log _{2}C_{\mu }$ is called the upper dimension of
$\Omega .$ We will assume that $\mu (\left\{ x\right\} )=0$ for all
$x\in \Omega ,$ (see Section \ref{sec2} below).

Our main objective will be to prove the equivalence between the
following non-collapsing condition: there exists a constant $b>0$
(called the non-collapsing constant) such that,
\begin{equation}
\inf_{x\in \Omega }\mu (B(x,1))=b,  \label{noncolap}
\end{equation}
and general versions of Haj{\l }asz-Besov type embeddings theorems\footnote{%
Precise definitions and properties concerning all the topics that
will appear in this Introduction and throughout the whole paper are
contained in Section \ref{section2} below.}. Since if $\mu (\Omega
)<\infty ,$ then obviously (\ref{noncolap}) holds, we will always
assume in what follows that $\mu (\Omega )=\infty .$

As in the Euclidean case, there are several equivalent ways to
define Haj{\l }asz-Besov spaces in the setting of doubling metric
measure spaces (see for example \cite{Amiran}, \cite{Amiran1},
\cite{HMY}, \cite{HIH}, \cite{MY}, \cite{Mu}, \cite{YZ} and the
references therein). In this paper, we will use the approach based
on modulus of smoothness for two reasons, (i) it just uses the
metric measure structure of $\Omega ,$ and (ii) its relation with
the real interpolation method.

We start, in order not to obscure the simplicity of the arguments,
by considering the version of our main result for Haj{\l }asz-Besov
whose modulus is generated by $L^{p}$-spaces, delaying the general
version to the last part of this introductory section. To this end,
we shall need some definitions.

Given $0<p<\infty $, the $L^{p}-$modulus of smoothness of a
measurable function $f$ in $\Omega $, is given by (see \cite{Amiran}
and \cite{CaoGri})
\begin{equation*}
\mathcal{E}_{p}(f,r):=\left( \int_{\Omega }\left( \frac{1}{\mu (B(x,r))}%
\int_{B(x,r)}\left| f(x)-f(y)\right| ^{p}d\mu (y)\right) d\mu
(x)\right) ^{1/p},\text{ }r>0,
\end{equation*}
and the homogeneous Haj{\l }asz-Besov space $\mathcal{\dot{B}}%
_{p,q}^{s}(\Omega ),0<s<1,$ $0<p<\infty $, $0<q\leq \infty ,$ is
defined as the set of measurable functions for which the Besov
seminorm
\begin{equation*}
\left\| f\right\| _{\mathcal{\dot{B}}_{p,q}^{s}(\Omega )}:=\left\{
\begin{array}{ll}
\left( \dint_{0}^{\infty }\left(
\frac{\mathcal{E}_{p}(f,t)}{t^{s}}\right)
^{q}\frac{dt}{t}\right) ^{1/q}, & 0<q<\infty , \\
\dsup_{t>0}t^{-s}\mathcal{E}_{p}(f,t), & q=\infty ,
\end{array}
\right.
\end{equation*}
is finite\footnote{%
In the Euclidean setting, applying Fubini's theorem it is easy to
verify that $\mathcal{E}_{p}(f,t)$ is equivalent to the classical
$L^{p}-$modulus of smoothness
\begin{equation*}
\omega _{p}(f,t)=\sup_{\left| h\right| \leq t}\left\|
f(x+h)-f(x)\right\| _{L^{p}}
\end{equation*}
(see \cite{Amiran}), therefore,
$\mathcal{\dot{B}}_{p,q}^{s}(\mathbb{R}^{n})$ coincides with the
classical definition.}.

For measurable functions $f:\Omega \rightarrow \mathbb{R},$ the
decreasing
rearrangement of $f$ is the decreasing function $f^{\ast }$ defined on $%
[0,\infty )$ by\footnote{%
See Section \ref{rearang} below.}
\begin{equation*}
f^{\ast }(t)=\inf \{s\geq 0:\mu _{f}(s)\leq t\},\text{ \ \ }t\geq 0.
\end{equation*}
where $\mu _{f}$ denotes the distribution function of $f.$

Associated to $f^{\ast },$ we consider the maximal function $f^{\ast
\ast
}(t)=\frac{1}{t}\int_{0}^{t}f^{\ast }(s)ds$, and for $0<\alpha \leq 1$ the $%
\alpha -$oscillation of $f^{\ast }$ defined by
\begin{equation*}
O(\left| f\right| ^{\alpha },t):=\left( \left| f\right| ^{\alpha
}\right) ^{\ast \ast }(t)-\left( \left| f\right| ^{\alpha }\right)
^{\ast }(t).
\end{equation*}

The $L^{p}$ version of our main result, is the following:

\begin{theorem}
\label{TeoLp}Let $\left( \Omega ,d,\mu \right) $ be doubling metric
measure space with doubling constant $C_{\mu }$ an upper dimension
$Q$. The following statements are equivalent

\begin{enumerate}
\item Condition (\ref{noncolap}) holds.

\item  Let $0<p\leq \infty ,\ 0<q\leq \infty $, $0<s<1$ and $\alpha =\min
(1,p),$ then there is a positive constant $C=C(C_{\mu },\alpha ,q)$
such that for all $f\in L^{\alpha }+L^{\infty },$ we have that
\begin{equation*}
\left( \int_{0}^{1}\left[ O(\left| f\right| ^{\alpha },t)^{\frac{1}{\alpha }%
}t^{\frac{1}{p}-\frac{s}{Q}}\right] ^{q}\frac{dt}{t}\right)
^{1/q}\leq C\left( \left\| f\right\|
_{\mathcal{\dot{B}}_{p,q}^{s}(\Omega )}+\left\| f\right\|
_{L^{\alpha }+L^{\infty }}\right) .
\end{equation*}
\end{enumerate}
\end{theorem}

\begin{remark}
\label{k11}In the setting of probability metric measure spaces the
technique to obtain pointwise estimates of the special differences
$O(|f|,t):=f^{\ast \ast }(t)-f^{\ast }(t),$ called the oscillation
of $f^{\ast }$, in terms of appropriate functionals depending on $f$
and the isoperimetric profile of the space has been developed by M.
Milman and J. Mart\'{i}n (see \cite{MM3}, \cite{MM4}, \cite{MM6})
and provide a considerable simplification in the theory of
embeddings of Sobolev spaces.
\end{remark}

We now turn to our main objective to present the corresponding
version of Theorem \ref{TeoLp} for Haj{\l }asz-Besov spaces defined
in a more abstract setting. The main idea is to change $L^{p}$ by
more general function spaces in the definition of the modulus of
smoothness, and the natural class seems to be the rearrangement
invariant quasi-Banach function spaces (\emph{q.r.i spaces} for
short) since, like Lebesgue spaces, q.r.i spaces have the property
equimeasurable functions $f,$ $g$ (i.e., such that $f^{\ast
}=g^{\ast }$) have equal quasi-norm\footnote{%
We shall refer to Banach function spaces that satisfies this
property as rearrangement invariant spaces (r.i. spaces for
short).}, (see Section \ref {Rea01} below).

Given a r.i. space $X$ on $\Omega $ and $0<r<\infty ,$ the $r-$%
convexification of $X$ is defined as $X^{(r)}=\{f:\left| f\right|
^{r}\in X\} $ endowed with the following quasi-norm $\left\|
f\right\| _{X^{(r)}}=\left\| \left| f\right| ^{r}\right\|
_{X}^{1/r}.$ It is plain that if $r\geq 1,$ then $X^{(r)}$ is a r.i.
space, however for $0<r<1$, the functional $\left\| \cdot \right\|
_{X^{(r)}}$ is not necessarily a norm. Notice that with this
definition we have that $(L^{p})^{(r)}=L^{pr}.$

Associated to $X^{(r)},$ we introduce a general variant of the
modulus of smoothness in the following way: Let $0<\alpha \leq 1$
and $f\in L^{\alpha
}(\Omega )+L^{\infty }(\Omega ).$ The $X^{(\alpha )}-$modulus of smoothness $%
E_{X^{(\alpha )}}:$ $X^{(\alpha )}\times (0,\infty )\rightarrow
(0,\infty )$ is defined by
\begin{equation*}
E_{X^{(\alpha )}}(f,r)=\left\| \nabla _{r}^{\alpha }f\right\|
_{X^{(\alpha )}},\text{\ }
\end{equation*}
where
\begin{equation*}
\nabla _{r}^{\alpha }f(x)=\left( \frac{1}{\mu
(B(x,r))}\int_{B(x,r)}\left| f(x)-f(y)\right| ^{\alpha }d\mu
(y)\right) ^{1/\alpha },\text{ \ \ }r>0.
\end{equation*}
The corresponding homogeneous Haj\l asz-Besov space\textbf{\ }$\dot{B}%
_{X^{(\alpha )},q}^{s}(\Omega ),$ $0<q\leq \infty $ and $0<s<1,$
defined by this modulus, consists of all functions $f\in L^{\alpha
}(\Omega )+L^{\infty }(\Omega )$ for which
\begin{equation*}
\left\| f\right\| _{\dot{B}_{X^{(\alpha )},q}^{s}(\Omega )}=\left(
\int_{0}^{\infty }\left( r^{-s}E_{X^{(\alpha )}}(f,r)\right) ^{q}\frac{dr}{r}%
\right) ^{1/q}
\end{equation*}
is finite (with the usual modification when $q=\infty $).

\begin{remark}
\label{modp}If $0<p\leq 1,$ then $E_{\left( L^{1}\right) ^{(p)}}(f,r)=%
\mathcal{E}_{p}(f,r),$ but if $1<p<\infty $, then by H\"{o}lder's
inequality, we get
\begin{align*}
E_{L^{p}}(f,r)& =\left( \int_{\Omega }\left( \frac{1}{\mu (B(x,r))}%
\int_{B(x,r)}\left| f(x)-f(y)\right| d\mu (y)\right) ^{p}d\mu
(x)\right)
^{1/p} \\
& \leq \left( \int_{\Omega }\left( \frac{1}{\mu
(B(x,r))}\int_{B(x,r)}\left|
f(x)-f(y)\right| ^{p}d\mu (y)\right) d\mu (x)\right) ^{1/p} \\
& =\mathcal{E}_{p}(f,r).
\end{align*}
However, we will see in Remark \ref{equal} below that both modulus
of smoothness produce the same Haj\l asz-Besov spaces.
\end{remark}

Our main result is the following:

\begin{theorem}
\label{k1}Let $\left( \Omega ,d,\mu \right) $ be doubling metric
measure space with doubling constant $C_{\mu }$ an upper dimension
$Q$. The following statements are equivalent

\begin{enumerate}
\item  Condition (\ref{noncolap}) holds.

\item Let $0<\alpha \leq 1,\ 0<q\leq \infty $ and $0<s<1.$ Then there is a
positive constant $C=C(C_{\mu },\alpha ,q)$ such that for any r.i.
space $X$ on $\Omega $ we have that
\begin{equation}
\left( \int_{0}^{1}\left( O(\left| f\right| ^{\alpha },t)^{1/\alpha }(t)%
\frac{\phi _{X^{(\alpha )}}(t)}{t^{s/Q}}\right)
^{q}\frac{dt}{t}\right) ^{1/q}\leq C\left( \left\| f\right\|
_{\dot{B}_{X^{(\alpha )},q}^{s}}+\left\| f\right\| _{L^{\alpha
}+L^{\infty }}\right) , \label{ineoscil}
\end{equation}
where $\phi _{X^{(\alpha )}}$ denotes the fundamental function of $%
X^{(\alpha )}$ (see Section \ref{Rea01} below).
\end{enumerate}
\end{theorem}

\begin{remark}
Theorem \ref{TeoLp} follows from Remark \ref{modp} and Theorem
\ref{k1}
considering $X=L^{p}$ and $\alpha =1$ if $p\geq 1,$ or $X=L^{1/p}$ and $%
\alpha =p$ if $p<1.$
\end{remark}

\bigskip

Theorem \ref{k1} represents an improvement over the known results
for two reasons, first does not require any Poincar\'{e} inequality
assumption or
topological requirement such as being uniformly perfect or being and $RD-$%
space, etc. (see for example \cite{AGH}, \cite{Amiran}, \cite{MaOr},
\cite {MY} or \cite{YZ} ) secondly it can be applied, following the
philosophy used in \cite{Curbera}, to q.r.i. spaces. Let us explain
this second point in more detail.

Given a r.i. space $X,$ the scale $X^{(\alpha )}$, $0<\alpha \leq
1$, has quasi-Banach spaces, but $\left( X^{(\alpha )}\right)
^{(1/\alpha )}=X$ is Banach. This observation is crucial if we want
to apply Theorem \ref{k1} to more general scales. The main idea is
to require that if $X$ is a q.r.i. space, then on the scale
$\{X^{(r)}\}_{1\leq r<\infty }$ there may be quasi-Banach spaces,
but there are also Banach spaces. This fact can be translated into a
convexity assumption about the space in the following way:

A q.r.i. space $X$ for which there is $0<\alpha \leq 1,$ such that
the functional $\left\| \cdot \right\| _{X^{(1/\alpha )}}$ is
equivalent to a norm is called an $\alpha $\emph{-convex space} (see
Section \ref{Reapcon} below). There are examples of q.r.i. spaces
that are not $\alpha -$convex for any $0<\alpha \leq 1$ (see
\cite{JS}) but they are very rare as Grafakos and Kalton said (see
\cite{GK}) ''all practical quasi-Banach rearrangement invariant
spaces are $\alpha -$convex for some $0<\alpha \leq 1",$ for this
reason throughout the paper we will restrict ourselves to work with $\alpha $%
-convex spaces.

The corresponding version of Theorem \ref{k1} for $\alpha $-convex
spaces reads as follows.

\begin{theorem}
\label{kk1}Let $\left( \Omega ,d,\mu \right) $ be doubling metric
measure space with doubling constant $C_{\mu }$ an upper dimension
$Q$. Let $X$ be an $\alpha -$convex q.r.i space ($0<\alpha \leq 1)$.
The following statements are equivalent.

\begin{enumerate}
\item The non-collapsing condition holds.

\item  Let $0<q\leq \infty $ and $0<s<1.$ Then there is a positive constant $%
C=C(C_{\mu },\alpha ,q)$ such that
\begin{equation}
\left( \int_{0}^{1}\left( O(\left| f\right| ^{\alpha },t)^{1/\alpha }(t)%
\frac{\phi _{X}(t)}{t^{s/Q}}\right) ^{q}\frac{dt}{t}\right)
^{1/q}\leq C\left( \left\| f\right\| _{\dot{B}_{X,q}^{s}}+\left\|
f\right\| _{L^{\alpha }+L^{\infty }}\right) .  \label{ineoscil1}
\end{equation}
\end{enumerate}
\end{theorem}

\begin{remark}
The formulation given in Theorem \ref{k1} and the equivalent one
presented in the previous Theorem reflect that there are two
different points of view: suppose that we want to get estimates in
$L^{1,\infty }.$ The first formulation consists of looking at the
q.r.i space $X=L^{1,\infty }$ which has the property that for any
$0<\alpha <1,$ $X^{(1/\alpha )}=$ $L^{1/\alpha ,\infty }$ is a r.i.
space. This convexity allows us to apply Theorem \ref {kk1} to $X$.
Alternatively we can start from $X=L^{1/\alpha ,\infty }$
which is a r.i. space and by the first formulation get estimates in $%
X^{\left( p\right) }$ for all $0<p<1$, and in particular in
$X^{(\alpha )}=L^{1,\infty }.$
\end{remark}

The paper is organized as follows. In Section \ref{section2}, we
introduce the notation and the standard assumptions used in the
paper, we give the basic definitions and results we will need from
the theory or q.r.i. spaces, Haj\l asz-Sobolev and Haz\l asz-Besov
spaces. Section \ref{secinterpol} is devoted to interpolation, our
main result is Theorem \ref{TeoInterpol} where show that Haj\l
asz-Besov spaces $\dot{B}_{X^{(\alpha )},q}^{s}\left( \Omega \right)
$ are interpolation spaces between $X^{(\alpha )}$ and Haj\l
asz-Sobolev spaces $\dot{M}^{1,X^{(\alpha )}}$ (see definitions on
Section \ref{HazSobo} below), that is,
\begin{equation*}
\dot{B}_{X^{(\alpha )},q}^{s}(\Omega )=(X^{(\alpha
)},M^{1,X^{(\alpha )}})_{s,q}.
\end{equation*}
for $0<s<1$, $0<p<\infty $ and $0<q\leq \infty .$ This result is an
extension of \cite[Theorem 4.1]{HIH} where it was proved that for
classical
Haj\l asz-Sobolev spaces $M^{1,p}$ one has that\footnote{%
In case of $p>1$, $q>1$ this was earlier obtained in \cite{Amiran}
under some additional assumptions.}
\begin{equation*}
\mathcal{\dot{B}}_{p,q}^{s}(\Omega )=\left(
L^{p},\dot{M}^{1,p}\right) _{s,q}.
\end{equation*}
Section \ref{section4} contains the proofs of our main results,
Theorems \ref {k1} and \ref{kk1}.

Finally, in Section \ref{sec05} we will derive embeddings results
from inequalities (\ref{ineoscil}) or (\ref{ineoscil1}) our main
difficulty will be that the left hand side of these inequalities
depends neither on the growth of $(|f|^{\alpha })^{\ast }$ nor on
$(|f|^{\alpha })^{\ast \ast }$ but rather on the oscillation, this
will solved by showing (see Lemma \ref {pesos}) that there is a
weight $w$ (which depends on $q$ and $\alpha )$ such that for all
$f\in L^{\alpha }+L^{\infty }$ we get
\begin{equation*}
\left( \int_{0}^{1}\left( f^{\ast }(t)w(t)\right)
^{q}\frac{dt}{t}\right) ^{1/q}\leq C\left( \left( \int_{0}^{1}\left(
O(\left| f\right| ^{\alpha
},t)^{1/\alpha }(t)\frac{\phi _{X}(t)}{t^{s/Q}}\right) ^{q}\frac{dt}{t}%
\right) ^{1/q}+\left\| f\right\| _{L^{\alpha }+L^{\infty }}\right) .
\end{equation*}
We also analyze the role that the upper dimension\footnote{%
Sobolev-Besov embedding for $\dot{B}_{p,q}^{s}(\mathbb{R}^{n})$
spaces has a different character if $0<p<\frac{s}{n}$,
$p=\frac{s}{n}$ or $p>\frac{s}{n}.$ In the doubling metric context
the counterpart of the dimension $n$ is given by the upper dimension
$Q$.} $Q$ and the parameter $\alpha $ play. Let us to explain this
second point in more detail. Consider the scale of Lebesgue
spaces $\left\{ L^{p}\right\} _{0<p<\infty },$ taking into account that $%
L^{p}$ is $\min (1,p)-$convex and that $\phi _{L^{p}}(t)=t^{1/p},$
(\ref {ineoscil1}) reads
\begin{equation*}
\left( \int_{0}^{1}\left( O(\left| f\right| ^{\min
(1,p)},t)^{\frac{1}{\min
(1,p)}}(t)t^{\frac{1}{p}-\frac{s}{Q}}\right) ^{q}\frac{dt}{t}\right)
^{1/q}\leq C\left( \left\| f\right\| _{\dot{B}_{p,q}^{s}}+\left\|
f\right\| _{L^{\min (1,p)}+L^{\infty }}\right) .
\end{equation*}
This inequality implies (see Theorem \ref{Teolorentzlog} below)

\begin{enumerate}
\item  If $s<\frac{Q}{p},$ then
\begin{equation*}
\left( \int_{0}^{1}\left[ t^{\frac{1}{p}-\frac{s}{Q}}f^{\ast }(t)\right] ^{q}%
\frac{dt}{t}\right) ^{1/q}\leq C\left( \left\| f\right\| _{\dot{B}%
_{p,q}^{s}}+\left\| f\right\| _{L^{\min (1,p)}+L^{\infty }}\right) .
\end{equation*}

\item  If $s=\frac{Q}{p}$ and $\min (1,p)<q<\infty ,$ then
\begin{equation}
\left( \int_{0}^{1}\left[ \frac{f^{\ast }(t)}{\left( 1+\ln \frac{1}{t}%
\right) ^{\frac{1}{\min (1,p)}}}\right] ^{q}\frac{dt}{t}\right)
^{1/q}C\left( \left\| f\right\| _{\dot{B}_{p,q}^{s}}+\left\|
f\right\| _{L^{\min (1,p)}+L^{\infty }}\right) .  \label{z1}
\end{equation}

\item  If $s=\frac{Q}{p},$ and $q\leq \min (1,p)$ or $s>\frac{Q}{p},$ then
\begin{equation*}
\left\| f\right\| _{L^{\infty }}\leq C\left( \left\| f\right\| _{\dot{B}%
_{p,q}^{s}}+\left\| f\right\| _{L^{\min (1,p)}+L^{\infty }}\right) .
\end{equation*}
\end{enumerate}

If instead of $\left\{ L^{p}\right\} _{0<p<\infty },$ we consider
the scale of weak-Lebesgue spaces $\left\{ L^{p,\infty }\right\}
_{0<p<\infty }$, then
$\phi _{L^{p,\infty }}(t)=t^{1/p},$ but for $0<p\leq 1\ $the space $%
L^{p,\infty }$ is $\alpha -$convex only if $\alpha <p$, thus
inequality (\ref {ineoscil1}) is now
\begin{equation*}
\left( \int_{0}^{1}\left( O(\left| f\right| ^{\alpha },t)^{1/\alpha }(t)t^{%
\frac{1}{p}-\frac{s}{Q}}\right) ^{q}\frac{dt}{t}\right) ^{1/q}\leq
C\left( \left\| f\right\| _{\dot{B}_{X,q}^{s}}+\left\| f\right\|
_{L^{\alpha }+L^{\infty }}\right) .
\end{equation*}
From here we obtain\footnote{%
As far we know this part is new even in the Euclidean case.} (see
Theorem \ref{Teolorentzlog} below)

\begin{enumerate}
\item  If $s<\frac{Q}{p},$ then
\begin{equation*}
\left( \int_{0}^{1}\left[ t^{\frac{1}{p}-\frac{s}{Q}}f^{\ast }(t)\right] ^{q}%
\frac{dt}{t}\right) ^{1/q}\leq C\left( \left\| f\right\| _{\dot{B}%
_{L^{p,\infty },q}^{s}}+\left\| f\right\| _{L^{\min (1,p)}+L^{\infty
}}\right) .
\end{equation*}

\item  If $s=\frac{Q}{p},$ then

\begin{enumerate}
\item  If $1<p,$ then
\begin{equation*}
\left( \int_{0}^{1}\left[ \frac{f^{\ast }(t)}{\left( 1+\ln \frac{1}{t}%
\right) }\right] ^{q}\frac{dt}{t}\right) ^{1/q}\leq C\left( \left\|
f\right\| _{\dot{B}_{L^{p,\infty },q}^{s}}+\left\| f\right\|
_{L^{\min (1,p)}+L^{\infty }}\right) .
\end{equation*}

\item  If $0<p\leq 1,$ then for any $0<\alpha <p$, we have that
\begin{equation}
\left( \int_{0}^{1}\left[ \frac{f^{\ast }(t)}{\left( 1+\ln \frac{1}{t}%
\right) ^{\frac{1}{\alpha }}}\right] ^{q}\frac{dt}{t}\right)
^{1/q}\leq C_{\alpha }\left( \left\| f\right\|
_{\dot{B}_{L^{p,\infty },q}^{s}}+\left\| f\right\| _{L^{\alpha
}+L^{\infty }}\right) .  \label{z2}
\end{equation}
\end{enumerate}

\item  If $s=\frac{Q}{p},$ and $q\leq \min (1,p)$ or $s>\frac{Q}{p},$ then
\begin{equation*}
\left\| f\right\| _{L^{\infty }}\leq C\left( \left\| f\right\| _{\dot{B}%
_{L^{p,\infty },q}^{s}}+\left\| f\right\| _{L^{\min (1,p)}+L^{\infty
}}\right) .
\end{equation*}
\end{enumerate}

\begin{remark}
Notice that $\phi _{L^{p}}(t)=\phi _{L^{p,\infty }}(t)$, however embeddings (%
\ref{z1}) and (\ref{z2}) have different expressions when $0<p\leq
1$, this is due to the fact that if $0<p\leq 1,$ then $\left(
L^{p}\right) ^{\left( 1/p\right) }=L^{1}$ is a Banach spaces but
$\left( L^{p,\infty }\right) ^{\left( 1/p\right) }=L^{1,\infty }$
not.
\end{remark}

\begin{remark}
The embedding of Besov spaces $B_{p,q}^{s}(\Omega )$, $p>1,q\geq 1$,
into Lorentz spaces was proved in \cite{Amiran} under the assumption
that $\Omega $ is $Q-$Ahlfors regular, i.e. if there is a constant
$c_{Q}$ $\geq 1$ such that
\begin{equation*}
\frac{r^{Q}}{c_{Q}}\leq \mu (B(x,r))\leq c_{Q}r^{Q}
\end{equation*}
and supports a $(1,p)-$Poincar\'{e} inequality. For $Q-$Ahlfors
regular spaces and $0<p<\infty ,\ 0<q\leq \infty $ was obtained in
\cite{HIH}.
\end{remark}

\

Throughout the paper, we denote by $C$ a positive constant which is
independent of the main parameters, but which may vary from line to
line.
The symbol $f\preceq g$ means that $f\leq c g$ for some $c > 0$. If $%
f\preceq g$ and $g\preceq f$ we then write $f\simeq g$.

\section{Notation and preliminaries\label{section2}}

In this section we establish some further notation and background
information and we provide more details about metrics spaces and
function spaces in which we will be working with.

\subsection{Metric spaces\label{sec2}}

Let $(\Omega ,d)$ be a metric space. As usual a ball $B$ in $\Omega
$ with a center $x$ and radius $r>0$ is a set $B=B(x,r):=\{y\in
\Omega ;d(x,y)<r\}$. Throughout the paper by a metric measure space
we mean a triple $(\Omega ,d,\mu )$, where $\mu $ is a Borel measure
on $(\Omega ,d)$ such $0<\mu (B)<\infty $, for every ball $B$ in
$\Omega $, we also assume that $\mu (\Omega )=\infty $ and $\mu
(\left\{ x\right\} )=0$ for all $x\in \Omega .$

We say that $\left( \Omega ,d,\mu \right) $ is a \textbf{doubling
metric space}, if $\mu $ is a \textbf{doubling measure} on $\Omega
$, i.e. $\mu $ satisfies that there exists a constant $C>1$, such
that
\begin{equation}
0<\mu (B(x,2r))\leq C\mu (B(x,r))<\infty  \label{dobla}
\end{equation}
for every ball $B,$ for all $x\in \Omega $ and $r>0.$ The smallest constant $%
C$ for which (\ref{dobla}) is satisfied will be called the doubling
constant and will be denoted by $C_{\mu }.$ The upper dimension of
$\Omega $ defined by
\begin{equation*}
Q=\log _{2}C_{\mu }.
\end{equation*}
By means of an iteration of condition (\ref{dobla}) we get (see for
example \cite[Lemma 4.7]{Ha1} for the details)
\begin{equation}
\mu (B(x,r))\geq \left( \frac{r}{4R}\right) ^{Q}\mu (B(y,R))
\label{iteradobla}
\end{equation}
whenever $x,y\in \Omega $ satisfy $B(x,r)\subset B(y,R)$ and
$0<r\leq R<\infty .$

We say that $\left( \Omega ,d,\mu \right) $ satisfies the \textbf{%
non-collapsing condition} if there is a constant $b>0,$ called the
non-collapsing constant, such that
\begin{equation*}
\inf_{x\in \Omega }\mu (B(x,1))=b.
\end{equation*}

\subsection{Rearrangements of functions\label{rearang}}

Let $\left( \Omega ,d,\mu \right) $ be a doubling metric measure
space. For measurable functions $f:\Omega \rightarrow \mathbb{R},$
the distribution function of $f$ is given by
\begin{equation*}
\mu _{f}(t)=\mu \{x\in {\Omega }:\left| f(x)\right| >t\}\text{ \ \ \ \ }%
(t>0).
\end{equation*}
The \textbf{decreasing rearrangement} of $f$ is the function
decreasing function $f^{\ast }$ defined on $[0,\infty )$ by
\begin{equation*}
f_{\mu }^{\ast }(t)=\inf \{s\geq 0:\mu _{f}(s)\leq t\},\text{ \ \
}t\geq 0.
\end{equation*}
The main property if $f^{\ast }$ of that it is equimeasurable with
$f$, i.e.
\begin{equation*}
\mu \{x\in {\Omega }:\left| f(x)\right| >t\}=\left| \{s\in \lbrack
0,\infty ):f^{\ast }(s)>t\}\right| .
\end{equation*}
For every positive $\alpha $ we have $\left( \left| f\right|
^{\alpha }\right) ^{\ast }(s)=\left( f^{\ast }(s)\right) ^{\alpha
},$ and if $|g|\leq |f|$ $\mu -$almost everywhere on $\Omega $, then
$g^{\ast }\leq f^{\ast }$. A detailed treatment of rearrangements
may be found in Bennett and Sharpley \cite{BS}.

A basic property of rearrangements is the Hardy-Littlewood
inequality which tells us that, if $f$ and $g$ are two $\mu
$-measurable functions on $\Omega $, then
\begin{equation}
\int_{{\Omega }}|f(x)g(x)|\,d\mu \leq \int_{0}^{\infty }f^{\ast
}(t)g^{\ast }(t)\,dt.  \label{HL}
\end{equation}
in particular, for any $\mu -$measurable set $E\subset \Omega $
\begin{equation*}
\int_{E}\left| f(x)\right| d\mu \leq \int_{0}^{\mu (E)}f^{\ast
}(s)ds.
\end{equation*}
Since $f^{\ast }$ is decreasing, the maximal function $f^{\ast \ast }$ of $%
f^{\ast },$ defined by
\begin{equation*}
f^{\ast \ast }(t)=\frac{1}{t}\int_{0}^{t}f^{\ast }(s)ds,
\end{equation*}
is also decreasing, and
\begin{equation*}
f^{\ast }\leq f^{\ast \ast }.
\end{equation*}
Notice that
\begin{equation}
\frac{\partial }{\partial t}f^{\ast \ast }(t)=-\left( \frac{f^{\ast
\ast }(t)-f^{\ast }(t)}{t}\right) .  \label{der2est}
\end{equation}

We single out two subadditivity properties, if $f$ and $g$ are two $\mu -$%
measurable functions on $\Omega $, then for $t>0$

\begin{equation}
\left( f+g\right) ^{\ast }(2t)\leq f^{\ast }(t)+g^{\ast }(t)
\label{xxx}
\end{equation}
and

\begin{equation}
\left( f+g\right) ^{\ast \ast }(t)\leq f^{\ast \ast }(t)+g^{\ast
\ast }(t). \label{readoble}
\end{equation}

These facts will be used throughout this paper.

\subsection{Background on Rearrangement Invariant Spaces\label{Rea01}}

We recall briefly the basic definitions and conventions we use from
the theory of rearrangement invariant spaces and refer the reader to
\cite{BS} and \cite{KPS}, for a complete treatment. Let $({\Omega
},d,\mu )$ be a
metric measure space. Let $L^{0}(\Omega )$ be the space of all real-valued $%
\mu -$measurable functions on $\Omega $, with the topology of local
convergence in measure. A \textbf{quasi-Banach rearrangement
invariant function space} $X=X({\Omega })$ on ${\Omega }$ (a
\textbf{q.r.i space}) will be a subspace of $L^{0}$ equipped with a
quasi-norm $\left\| \cdot \right\| _{X}$ such that

\begin{enumerate}
\item  $X$ is complete (i.e. a quasi-Banach space) for $\left\| \cdot
\right\| _{X}.$

\item  The injection $X\rightarrow L^{0}$ is continuous.

\item  If $E$ is a set of finite measure, then $\chi _{E}\in X.$

\item  If $0\leq f_{n}\nearrow f$ $\mu -a.e.$ and $f\in X$ then $\left\|
f_{n}\right\| _{X}\nearrow \left\| f\right\| _{X}.$

\item  $\Vert f\Vert _{X}=\Vert g\Vert _{X}$ whenever $f^{\ast }=g^{\ast }.$
\end{enumerate}

If $\left\| \cdot \right\| _{X}$ is a norm, then we shall refer to
$X$ as a \textbf{r.i space.}

A r.i. space $X({\Omega })$ can be represented by an r.i. space on
the interval $(0,\infty ),$ with Lebesgue measure,
$\bar{X}=\bar{X}(0,\infty ),$ such that
\begin{equation*}
\Vert f\Vert _{X}=\Vert f^{\ast }\Vert _{\bar{X}},
\end{equation*}
for every $f\in X.$ A characterization of the norm $\Vert \cdot \Vert _{\bar{%
X}}$ is available (see \cite[Theorem 4.10 and subsequent
remarks]{BS}).

The associated space $X^{\prime }(\Omega )$ of a r.i. space
$X(\Omega )$ is the r.i. space of all measurable functions $h$ for
which the r.i. norm given by
\begin{equation}
\left\| h\right\| _{X^{\prime }(\Omega )}=\sup_{g\neq
0}\frac{\dint_{\Omega }\left| g(x)h(x)\right| d\mu }{\left\|
g\right\| _{X(\Omega )}}. \label{holhol}
\end{equation}
Note that by (\ref{holhol}), the generalized H\"{o}lder inequality

\begin{equation}
\int_{\Omega }\left| g(x)h(x)\right| d\mu \leq \left\| g\right\|
_{X(\Omega )}\left\| h\right\| _{X^{\prime }(\Omega )}
\label{holin}
\end{equation}
holds.

An important consequence of (\ref{HL}) and (\ref{holhol}) is the
Hardy-Littlewood-P\'{o}lya principle stating that
\begin{equation}
\int_{0}^{t}f^{\ast }(s)ds\leq \int_{0}^{t}g^{\ast }(s)ds\quad
\forall t>0\Rightarrow \left\| f\right\| _{X}\leq \left\| g\right\|
_{X}\text{ \ for any r.i. space }X.  \label{Hardy}
\end{equation}

A useful tool, in the study of a q.r.i. space $X$ is the
\textbf{fundamental function\ }of $X$ defined by
\begin{equation*}
\phi _{X}(t)=\left\| \chi _{E}\right\| _{X},
\end{equation*}
where $E$ is any measurable subset of $\Omega $ with $\mu (E)=t.$
This function is increasing with $\phi _{X}(0)=0.$ For example if
$X=L^{p},$ then $\phi _{L^{p}}(t)=t^{1/p}.$

In case that $X$ is a r.i. space, then $\phi _{X}$ is quasi-concave,
that is, $\phi _{X}(t)/t$ is decreasing. By renorming, if necessary,
we can always assume that $\phi _{X}$ is concave. Moreover,
\begin{equation}
\phi _{X^{\prime }}(s)\phi _{X}(s)=s.  \label{dual}
\end{equation}
\ Associated with an r.i. space $X$ there are some useful Lorentz
and Marcinkiewicz spaces, namely the Lorentz and Marcinkiewicz
spaces defined by the quasi-norms
\begin{equation*}
\left\| f\right\| _{M(X)}=\sup_{t}f^{\ast\ast}(t)\phi _{X}(t),\text{ \ \ }%
\left\| f\right\| _{\Lambda (X)}=\int_{0}^{\infty }f^{\ast }(t)d\phi
_{X}(t).
\end{equation*}
Notice that
\begin{equation*}
\phi _{M(X)}(t)=\phi _{\Lambda (X)}(t)=\phi _{X}(t),
\end{equation*}
and that
\begin{equation*}
\Lambda (X)\subset X\subset M(X).
\end{equation*}

\

Let $r>0$ and let $X$ be a q.r.i. space on $\Omega ;$ the $r-$\textbf{%
convexification} $X^{(r)}$ of $X,$ (see. \cite{LT} and \cite{JS} ) is defined%
\textbf{\ }by
\begin{equation*}
X^{(r)}=\{f:\left| f\right| ^{r}\in X\},\text{ \ \ }\left\|
f\right\| _{X^{(r)}}=\left\| \left| f\right| ^{r}\right\|
_{X}^{1/r}.
\end{equation*}

It follows also that $\phi _{X^{(r)}}(t)=\left( \phi _{X}(t)\right)
^{1/r}$. As we pointed in the introduction section, if $X$ is a r.i.
space and $r\geq 1$ then, $X^{(r)}$ still is a r.i. space but, in
general, for $0<r<1$, the space $X^{(r)}$ is not necessarily Banach.

\subsubsection{\ $\protect\alpha -$convex q.r.i. spaces\label{Reapcon}}

As mentioned, in the introduction section, we are interested in
q.r.i. spaces $X$ with the property that its $r-$convexification is
a Banach space for some large power $r.$ Here we will give some
equivalent characterizations for this class of spaces.

Let $X$ be a q.r.i space $X,$\ we will say that $X$ is $\alpha
-$convex for some $0<\alpha \leq 1,$ it there exists a positive
constant $C$ such that for all $f_{1},\cdots ,f_{n}\in X$ we have
\begin{equation}
\left\| \left( \sum_{j=1}^{n}\left| f_{j}\right| ^{\alpha }\right)
^{1/\alpha }\right\| _{X}\leq C\left( \sum_{j=1}^{n}\left\|
f_{j}\right\| _{X}^{\alpha }\right) ^{1/\alpha }.  \label{pconvex}
\end{equation}
For $0<q<\infty $ we say $X$ is $q-$concave if for some $C<\infty .$%
\begin{equation*}
\left( \sum_{j=1}^{n}\left\| f_{j}\right\| _{X}^{q}\right)
^{1/q}\leq C\left\| \left( \sum_{j=1}^{n}\left| f_{j}\right|
^{q}\right) ^{1/q}\right\| _{X}.
\end{equation*}
$X$ is called geometrically convex (see \cite{Kal1}) there exists a
positive constant $C$ such that for all $f_{1},\cdots ,f_{n}\in X$
we get
\begin{equation*}
\left\| \left| f_{1}\cdots f_{n}\right| ^{1/n}\right\| _{X}\leq
C\left( \prod_{j=1}^{n}\left\| f_{j}\right\| _{X}\right) ^{1/n}.
\end{equation*}
Finally, $X$ is said to be $L-$convex (see \cite{Kal}) if there is
$0<\varepsilon <1$ so that
if $0\leq f\in X$ with $\left\| f\right\| _{X}=1$ and $0\leq f_{i}\leq f,$ $%
i=1,\cdots ,n,$ satisfy
\begin{equation*}
\frac{1}{n}\sum_{i=1}^{n}f_{i}\geq (1-\varepsilon )f,
\end{equation*}
then
\begin{equation*}
\max_{1\leq i\leq n}\left\| f_{i}\right\| _{X}\geq \varepsilon .
\end{equation*}

\begin{proposition}
Let $X$ and $Y$ be q.r.i spaces.
\begin{enumerate}
\item[(i)] If $X$ $\alpha-$\emph{convex}, then the functional $\left\| \cdot \right\| _{X^{(1/\alpha)}}$ is equivalent to a norm.
\item[(ii)]
The following statements are equivalent:
\begin{enumerate}
 \item $X$ is $\alpha -$convex.
  \item $X$ is $L-$convex.
  \item $X$ is geometrically convex.
\end{enumerate}
\item[(iii)]  If $X$ is $q-$concave for some $0<q<\infty $, then $X$ is $L-$convex.
\item[(iv)]  If $X$ and $Y$ are $\alpha -$convex, then spaces $X+Y$ and $X\bigcap Y
$ are $\alpha -$convex q.r.i spaces.

\item[(v)] If $X$ is $\alpha -$convex, then
\begin{equation}
L^{\infty }(\Omega )\cap L^{\alpha }(\Omega )\subset X^{(\alpha
)}(\Omega )\subset L^{\alpha }(\Omega )+L^{\infty }(\Omega ),
\label{embbb}
\end{equation}
with continuous embedding.
\end{enumerate}
\end{proposition}

\begin{proof}
(i) If $f_{1},\cdots ,f_{n}\in X^{(1/\alpha )},$ then $g_{j}=\left|
f_{j}\right| ^{1/\alpha }\in X,$ thus from (\ref{pconvex}) it
follows
\begin{eqnarray*}
\left\| \sum_{j=1}^{n}\left| f_{j}\right| \right\| _{X^{(1/\alpha
)}} &=&\left\| \left( \sum_{j=1}^{n}\left| g_{j}\right| ^{\alpha
}\right)
^{1/\alpha }\right\| _{X}^{\alpha } \\
&\leq &C^{\alpha }\sum_{j=1}^{n}\left\| g_{j}\right\| _{X}^{\alpha } \\
&=&C^{\alpha }\sum_{j=1}^{n}\left\| f_{j}\right\| _{X^{1/\alpha}}.
\end{eqnarray*}
Now, it is easy to show that the functional
\begin{equation*}
\left\| f\right\| :=\inf \left\{ \sum_{j=1}^{n}\left\| f_{j}\right\|
_{X^{(1/\alpha )}}:\left| f\right| \leq \left( \sum_{j=1}^{n}\left|
f_{j}\right| ^{1/\alpha }\right) ^{\alpha }:f_{1},\cdots ,f_{n}\in
X^{(1/\alpha )}\right\}
\end{equation*}
is a r.i. norm on the space $X^{(1/\alpha )}$, which is equivalent
to the original quasi-norm.

(ii) Follows immediately because every $L-$convex q.r.i. space is $\alpha -$%
convex for some $0<\alpha \leq 1$ \cite[Theorem 2.2]{Kal}. On the
other hand, from \cite{Kal} and \cite{Kal1} we get that $X$ is
$L-$convex if and only if it is geometrically convex.

(iii) Was proved in \cite{Kal}, (iv) in \cite[Lemma 3.2]{AN} and (v)
in \cite {BHS}.
\end{proof}

It follows from the previous Proposition that given an $\alpha
-$convex q.r.i. space $X$, we have that $Y=X^{(1/\alpha )}$ is
equivalent to a r.i. space. Moreover, since obviously
\begin{equation*}
Y^{(\alpha )}=X
\end{equation*}
we get that all $\alpha -$convex q.r.i. spaces can be obtained as $\alpha -$%
convexification of r.i. spaces, for this reason
throughout the paper we will work with r.i. spaces or with $\alpha -$%
convexifications of r.i. spaces for some $0<\alpha \leq 1$.

\subsubsection{Examples}

We shall see below that all commonly arising examples of q.r.i
spaces (e.g. the $L^{p}$-spaces, Orlicz spaces, Lorentz spaces,
Lorentz-Zygmund spaces, generalized Lorentz-Zygmund spaces,
Marcinkiewicz spaces, etc.) are $\alpha - $convex.

\begin{enumerate}
\item  Classical Lorentz spaces: The spaces $L^{p,q}$ are defined by the
function quasi-norm
\begin{equation*}
\left\| f\right\| _{L^{p,q}}=\left( \int_{0}^{\infty }\left(
t^{1/p}f^{\ast }(t)\right) ^{q}\frac{dt}{t}\right) ^{1/q}
\end{equation*}
when $0<p,q<\infty ,$ and
\begin{equation*}
\left\| f\right\| _{L^{p,\infty }}=\sup_{0<t<\infty }t^{1/p}f^{\ast
}(t).
\end{equation*}

If $X=L^{p,q}$ , then $\left( L^{p,q}\right) ^{(r)}=L^{rp,rq}$ which
is equivalent to a Banach spaces for $r>\max \left\{ 1/p,1/q\right\}
.$

\item  Lorentz $\Lambda -$spaces: The Lorentz spaces $\Lambda ^{q}(w)$ are
defined by the functional
\begin{equation*}
\left\| f\right\| _{\Lambda ^{q}(w)}=\left( \int_{0}^{\infty
}f^{\ast }(t)^{q}w(t)dt\right) ^{1/q}
\end{equation*}
where $0<q<\infty $ and $w$ is a weight on $\left( 0,\infty \right) .$ If $%
w(t)=t^{q/p-1},$ then one obtain $\Lambda ^{q}(w)=L^{p,q}.$ If $%
w(t)=t^{q/p-1}\left( 1+\log ^{+}\frac{1}{t}\right) ^{q\beta },$ then $%
\Lambda ^{q}(w)=L^{p,q}(\log L)^{\beta }$ are the Lorentz-Zygmund
spaces
(see \cite{BS}) or if we take $w(t)=t^{q/p-1}\left( 1+\log ^{+}\frac{1}{t}%
\right) ^{q\beta }(1+\log ^{+}\log ^{+}\frac{1}{t})^{q\gamma }$ then $%
\Lambda ^{q}(w)=L^{p,q}(\log L)^{\beta }(\log \log L)^{\beta }$ are
the generalized Lorentz-Zygmund spaces.

If $X=\Lambda ^{q}(w)$ where the weight $w$ satisfies that there
exists large enough $r>1$ such that
\begin{equation*}
\int_{t}^{\infty }w(z)\frac{dz}{z^{r}}\leq \frac{C}{t^{r}}\int_{0}^{t}w(z)dz,%
\text{ \ \ }t>0,
\end{equation*}
then (see \cite{Sawyer}) $\Lambda ^{r}(w)$ is equivalent to a Banach
space$.$ This combined with the fact that $\left( \Lambda
^{1}(w)\right) ^{(p)}=\Lambda ^{p}(w),$ implies that $X^{(r)}$ is
equivalent to a r.i. spaces for large $r.$

\item  Marcinkiewicz spaces: Let $\varphi $ be an increasing concave
function with $\varphi (0^{+})=0.$ The Marcinkiewicz space
$M_{\varphi }$ is defined by the function norm
\begin{equation*}
\left\| f\right\| _{M_{\varphi }}=\sup_{t>0}\frac{\varphi (t)}{t}%
\int_{0}^{t}f^{\ast }(z)dz.
\end{equation*}
The space $M_{\varphi }$ is a r.i. space with fundamental function
$\varphi .
$ If $X$ is a r.i. space with fundamental function $\phi _{X},$ then $%
X\subset M_{\phi _{X}}.$ We also consider the q.r.i space $\tilde{M}%
_{\varphi }$ defined by the functional
\begin{equation*}
\left\| f\right\| _{\tilde{M}_{\varphi }}=\sup_{t>0}\varphi
(t)f^{\ast }(t).
\end{equation*}
Obviously $M_{\varphi }\subset \tilde{M}_{\varphi }.$ If
\begin{equation}
\frac{\varphi (t)}{t}\int_{0}^{t}\frac{1}{\varphi (z)}dz\leq C,
\label{marc}
\end{equation}
then $M_{\varphi }=\tilde{M}_{\varphi }.$ Moreover for any $r>1,$
$\left( \tilde{M}_{\varphi }\right) ^{(r)}$ is a r.i. space. Indeed,
direct
computation shows that $\left( \tilde{M}_{\varphi }\right) ^{(r)}=\tilde{M}%
_{\varphi ^{1/r}},$ and since $\varphi ^{1/r}$ is concave and
$\varphi ^{1/r}(t)/t^{1/r}$ is decreasing hence (\ref{marc}) holds.

For example if $X=L^{1},$ then $\phi _{X}=t$ and $M_{\phi _{X}}=L^{1}$ but $%
\tilde{M}_{\phi _{X}}=L^{1,\infty }$ that is not normable. However, $%
(L^{1,\infty })$ $^{(r)}$ is a Banach space for any $r>1.$

\item  Orlicz spaces: An Orlicz function is a continuous strict1y increasing
function $\Phi :[0,\infty )\rightarrow \lbrack 0,\infty )$ with
$\Phi (0)=0$ and satisfying the $\Delta _{2}-$condition$,$ i.e.
$\Phi (2x)\leq C\Phi (x).$

If we suppose (Matuszewska and Orlicz \cite{MO}) that for some $p>0$%
\begin{equation*}
\inf_{x,y\geq 1}\frac{\Phi (xy)}{\Phi (x)y^{p}}>0
\end{equation*}
then the Orlicz space $L^{\Phi }$ is defined by
\begin{equation*}
\left\| f\right\| _{L^{\Phi }}=\inf \left\{ \lambda >0:\int_{\Omega
}\Phi \left( \frac{\left| f(x)\right| }{\lambda }\right) d\mu
(x)\leq 1\right\} .
\end{equation*}
is a q.r.i. space.

$L^{\Phi }$ is $\alpha -$convex for $0<\alpha <\sup \left\{
p:\inf_{x,y\geq 1}\frac{\Phi (xy)}{\Phi (x)y^{p}}>0\right\} $ (see
\cite{LT}).
\end{enumerate}

\subsection{Haj{\l }asz-Sobolev and Haj{\l }asz-Besov spaces built on r.i.
spaces \label{HazSobo}}

Let $\left( \Omega ,d,\mu \right) $ be doubling metric measure
space. Let $X$ be a r.i. space on $\Omega \ $and $0<\alpha \leq 1.$

We say that a $\mu -$measurable function $f\in M^{1,X^{(\alpha
)}}\left( \Omega \right) $, if there exits a nonnegative measurable
function $g\in X^{(\alpha )}$ such that
\begin{equation}
\left| f(x)-f(y)\right| \leq d(x,y)\left( g(x)+g(y)\right) \text{ \
\ }\mu -a.e.\text{ \ }x,y\in \Omega .  \label{sgradi}
\end{equation}
A function $g$ satisfying (\ref{sgradi}) will be called a $1-$gradient of $f$%
. We denote by $D(f)$ the collection of all $1-$gradients of $f$.
The \textbf{homogeneous Haj\l asz Sobolev} space
$\dot{M}^{1,X^{(\alpha )}}(\Omega )$ consists of measurable
functions $f\in L^{\alpha }+L^{\infty }$ for which

\begin{equation*}
\left\| f\right\| _{\dot{M}^{1,X^{(\alpha )}}\left( \Omega \right)
}=\inf_{g\in D(f)}\left\| g\right\| _{X^{(\alpha )}}
\end{equation*}
is finite.

The \textbf{Haj\l asz-Sobolev} space $M^{1,X^{(\alpha )}}(\Omega )$ is $\dot{%
M}^{1,X^{(\alpha )}}(S)\cap X^{(\alpha )}$ equipped with the
quasi-norm
\begin{equation*}
\left\| f\right\| _{M^{1,X^{(\alpha )}}\left( \Omega \right)
}=\left\| f\right\| _{X^{(\alpha )}}+\left\| f\right\|
_{\dot{M}^{1,X^{(\alpha )}}\left( \Omega \right) }.
\end{equation*}
When $\alpha =1$ we will write $\dot{M}^{1,X}$ (resp. $M^{1,X}).$

\begin{remark}
The definition formulated above is motivated by the Haj\l asz
approach to the definition of Sobolev spaces on a metric measure
space (see \cite{Ha2} and \cite{Ha1}) where $M^{1,p}(\Omega )$ was
defined as the set of measurable functions $f$ for which
\begin{equation*}
\left\| f\right\| _{\dot{M}^{1,p}\left( \Omega \right) }=\inf_{g\in
D(f)}\left\| g\right\| _{L^{p}}
\end{equation*}
is finite\footnote{%
For $p>1$, $M^{1,p}(\mathbb{R}^{n})=W^{1,p}(\mathbb{R}^{n})$
\cite{Ha2}, \cite{Ha1}, whereas for $n/(n+1)<p\leq 1$,
$M^{1,p}(\mathbb{R}^{n})$ coincides with the Hardy-Sobolev space
$H^{1,p}(\mathbb{R}^{n})$ \cite[Theorem 1]{KoSak}.}. Based on this
definition spaces $M^{1,X^{(r)}}$ appear naturally when replacing
the Lebesgue norm by the quasi-norm $\left\| \cdot \right\|
_{X^{(r)}}.$

This generalization have been previously considered in some
particular cases, for example Tuominen \cite{Tuo} considered the
Orlicz case. Recently Heikkinen and Karak in \cite{HKa} have studied
in detail Orlicz-Has\l asz-Sobolev spaces. Costea and Miranda
\cite{CosMi} studied the Lorentz case. Mal\'{y}, in \cite{Mal2} and
\cite{Mal1}, investigated spaces associated with a general
quasi-Banach function lattice. Spaces $M^{1,Z}\left( \Omega \right)
$ where $Z$ is a r.i. space was considered in \cite{MaOr1}.
\end{remark}

\

Given $f\in L^{\alpha }(\Omega )+L^{\infty }(\Omega )$, we define the $%
X^{(\alpha )}-$modulus of continuity $E_{X^{(\alpha )}}:$
$X^{(\alpha )}\times (0,\infty )\rightarrow (0,\infty )$ by\ \
\begin{equation*}
E_{X^{(\alpha )}}(f,r)=\left\| \nabla _{r}^{\alpha }f\right\|
_{X^{(\alpha )}},\text{\ }
\end{equation*}
where
\begin{equation*}
\nabla _{r}^{\alpha }f(x)=\left( \frac{1}{\mu
(B(x,r))}\int_{B(x,r)}\left| f(x)-f(y)\right| ^{\alpha }d\mu
(y)\right) ^{1/\alpha },\text{ \ \ }r>0.
\end{equation*}
Let $0<q\leq \infty $ and $0<s<1.$ The \textbf{homogeneous Haj\l asz-Besov
space} $\dot{B}_{X^{(\alpha )},q}^{s}(\Omega )$ consists of all functions $%
f\in L^{\alpha }(\Omega )+L^{\infty }(\Omega )$ for which
\begin{equation*}
\left\| f\right\| _{\dot{B}_{X^{(\alpha )},q}^{s}(\Omega )}=\left(
\int_{0}^{\infty }\left( r^{-s}E_{X^{(\alpha )}}(f,r)\right) ^{q}\frac{dr}{r}%
\right) ^{1/q}
\end{equation*}
is finite (with the usual modification when $q=$ $\infty $).

Similarly the \textbf{Haj\l asz-Besov space} $B_{X^{(\alpha
)},q}^{s}(\Omega ) $ is $\dot{B}_{X^{(\alpha )},q}^{s}(\Omega )\cap
X^{(\alpha )}$ equipped with the quasi-norm
\begin{equation*}
\left\| f\right\| _{B_{X^{(\alpha )},q}^{s}}=\left\| f\right\|
_{X^{(\alpha )}}+\left\| f\right\| _{\dot{B}_{X^{(\alpha
)},q}^{s}\left( \Omega \right) }.
\end{equation*}
When $\alpha =1$ we write $\dot{B}_{X,q}^{s}$ (resp. $B_{X,q}^{s}).$

\section{Interpolation theorems for {Haj\l asz-Besov} spaces\label%
{secinterpol}}

In this section, we prove new interpolation for Haj\l asz-Besov
spaces. Let us start recalling some essential definitions and
properties of the real interpolation theory; see, for example, the
classical references \cite{BS}, \cite{BL} for the details.

Let $X_{0}$ and $X_{1}$ be (quasi-semi)normed spaces continuously
embedded into a topological vector space $\mathcal{X}$. For every
$f\in X_{0}+X_{1}$ and $t>0$, the $K-$functional is

\begin{equation*}
K(f,t;X_{0},X_{1})=\inf_{f=f_{0}+f_{1}}\left\{ \left\| f_{0}\right\|
_{X_{0}}+t\left\| f_{1}\right\| _{X_{1}}\right\} .
\end{equation*}

Let $0<s<1$ and $0<q\leq \infty $. The interpolation space $%
(X_{0},X_{1})_{s,q}$ consists of functions $f\in X_{0}+X_{1}$ for
which
\begin{equation*}
\left\| f\right\| _{(X_{0},X_{1})_{s,q}}=
\begin{cases}
\left( \int_{0}^{\infty }\left( t^{-s}K(f,t;X_{0},X_{1})\right) ^{q}\frac{dt%
}{t}\right) ^{1/q}, & \text{if }0<q<\infty , \\
\sup_{t}t^{-s}K(f,t;X_{0},X_{1}), & \text{if }q=\infty ,
\end{cases}
\end{equation*}
is finite.

Our main result in this section gives us the relation between
$K-$functional
for the couple $(X^{(\alpha )},\dot{M}^{1,X^{(\alpha )}})$ and the $%
X^{(\alpha )}-$modulus of smoothness.

\begin{theorem}
\label{TeoInterpol}Let $\left( \Omega ,d,\mu \right) $ be a doubling
measure metric space. Let $X$ be an r.i. space and $0<\alpha \leq
1,$ then that for all $f\in L^{\alpha }+L^{\infty \text{ }}$and
$t>0$,
\begin{equation}
E_{X^{(\alpha )}}(f,t)\preceq K(f,t,X^{(\alpha
)},\dot{M}^{1,X^{(\alpha )}})\preceq \left( \sum_{j=0}^{\infty
}2^{-j\alpha }E_{X^{(\alpha )}}(f,2^{j}t)^{\alpha }\right)
^{1/\alpha }  \label{kf1}
\end{equation}
and
\begin{equation}
K(f,t,X^{(\alpha )},M^{1,X^{(\alpha )}})\simeq K(f,t,X^{(\alpha )},\dot{M}%
^{1,X^{(\alpha )}})+\min (1,t)\left\| f\right\| _{X^{(\alpha )}}.
\label{kf2}
\end{equation}
\end{theorem}

This result will achieved through an appropriate modification of the
proof given in {\cite[Theorem 4.1]{HIH}} for the couple
$(L^{p},\dot{M}^{1,p})$, where the authors obtained a similar
expression for (\ref{kf1}) considering the modulus $\mathcal{E}_{p}$
instead of $E_{X^{(\alpha )}}$. The key point in our approach will
be the next result.

\begin{lemma}
\label{Lfubi}Let $X$ be an r.i. space and $0<\alpha \leq 1.$ Given
$f\in L^{\alpha }+L^{\infty }$ and $r>0$ we define
\begin{equation*}
T_{r}^{\alpha }f(x)=\left( \frac{1}{\mu (B(x,r))}\int_{B(x,r)}\left|
f(y)\right| ^{\alpha }d\mu (y)\right) ^{1/\alpha }.
\end{equation*}
Then the family of operators $\left\{ T_{r}^{\alpha }\right\}
_{r>0}$ is uniformly bounded on $X^{(\alpha )}.$
\end{lemma}

\begin{proof}
We claim that the family of operators $\left\{ T_{r}^{\alpha
}\right\} _{r>0}\ $is uniformly bounded on $L^{\infty \text{ }}$and
on $L^{\alpha }.$ Obviously
\begin{equation}
\left\| T_{r}^{\alpha }f(x)\right\| _{L^{\infty \text{ }}}\leq
\left\| f\right\| _{L^{\infty \text{ }}}.  \label{a1}
\end{equation}
On the other hand, by Fubini's theorem, we get
\begin{eqnarray*}
\left\| T_{r}^{\alpha }f\right\| _{L^{\alpha }}^{\alpha }
&=&\int_{\Omega }\left( \frac{1}{\mu (B(x,r))}\int_{B(x,r)}\left|
f(y)\right| ^{\alpha }d\mu
(y)\right) d\mu (x) \\
&=&\int_{\Omega }\left| f(y)\right| ^{\alpha }\left( \int_{B(y,r)}\frac{1}{%
\mu (B(x,r))}d\mu (x)\right) d\mu (y).
\end{eqnarray*}
Using the doubling property and the fact that $B(x,r)\subset
B(y,2r)$ whenever $y\in B(x,r)$, we conclude that
\begin{eqnarray*}
\int_{B(y,r)}\frac{1}{\mu (B(x,r))}d\mu (x) &\leq &C_{\mu }\int_{B(y,r)}%
\frac{1}{\mu (B(x,2r))}d\mu (x) \\
&\leq &C_{\mu }\int_{B(y,r)}\frac{1}{\mu (B(y,r))}d\mu (x) \\
&=&C_{\mu }.
\end{eqnarray*}
Thus, for all $r>0$
\begin{equation}
\left\| T_{r}^{\alpha }f\right\| _{L^{\alpha }}^{\alpha }\leq C_{\mu
}\int_{\Omega }\left| f(y)\right| ^{\alpha }d\mu (y)=C_{\mu }\left\|
f\right\| _{L^{\alpha }}^{\alpha }.  \label{a2}
\end{equation}
By combining (\ref{a1}) and (\ref{a2}) with the definition of the
$K$ functional, we obtain
\begin{equation*}
K(T_{r}^{\alpha }f,t^{1/\alpha },L^{\alpha },L^{\infty })\preceq
K(f,t^{1/\alpha },L^{\alpha },L^{\infty }),\text{ \ \ }(t>0).
\end{equation*}
Since (see \cite[Theorem 5.1]{BL})
\begin{equation*}
K(g,t^{1/\alpha },L^{\alpha },L^{\infty })\simeq \left(
\int_{0}^{t}\left( g^{\ast }(s)\right) ^{\alpha }ds\right)
^{1/\alpha }
\end{equation*}
we have
\begin{equation*}
\int_{0}^{t}\left( \left( T_{r}^{\alpha }f\right) ^{\ast }(s)\right)
^{\alpha }ds\preceq \int_{0}^{t}\left( f^{\ast }(s)\right) ^{\alpha
}ds
\end{equation*}
which by (\ref{Hardy}) implies that
\begin{equation*}
\left\| T_{r}^{\alpha }f\right\| _{X}\preceq \left\| \left| f\right|
^{\alpha }\right\| _{X},
\end{equation*}
as we was wished to show.
\end{proof}

\begin{proof}
(\textbf{of Theorem \ref{TeoInterpol})} We begin with the first
inequality
of (\ref{kf1}). Let $f=g+h$, where $g\in X^{(\alpha )}$, $h\in \dot{M}%
^{1,X^{(\alpha )}}$ and let $t>0$. Taking into account that
$0<\alpha \leq 1, $ we get
\begin{equation*}
\left| g(x)-g(y)\right| ^{\alpha }\leq \left| g(x)\right| ^{\alpha
}+\left| g(y)\right| ^{\alpha }
\end{equation*}
thus
\begin{equation*}
\left( \frac{1}{\mu (B(x,t))}\int_{B(x,t)}\left| g(x)-g(y)\right|
^{\alpha }d\mu (y)\right) ^{1/\alpha }\preceq \left| g(x)\right|
+T_{r}^{\alpha }g(x),\
\end{equation*}
consequently, by Lemma \ref{Lfubi},
\begin{equation*}
E_{X^{(\alpha )}}(g,t)\preceq \left\| g\right\| _{X^{(\alpha
)}}+\left\| T_{r}^{\alpha }g\right\| _{X^{(\alpha )}}\preceq \left\|
g\right\| _{X^{(\alpha )}}.
\end{equation*}
On the other hand, since $h\in \dot{M}^{1,X^{(\alpha )}},$ by the
definition of the $1-$gradient, if $\varrho \in D(h)\cap X^{(\alpha
)},$ then
\begin{align*}
\nabla _{r}^{\alpha }h(x)& =\left( \frac{1}{\mu (B(x,r))}%
\int_{B(x,r)}|h(x)-h(y)|^{\alpha }d\mu (y)\right) ^{1/\alpha } \\
& \leq \left( \frac{1}{\mu (B(x,t))}\int_{B(x,t)}d(x,y)^{\alpha
}|\varrho
(x)+\varrho (y)|^{\alpha }d\mu (y)\right) ^{1/\alpha } \\
& \leq \left( \frac{1}{\mu (B(x,t))}\int_{B(x,t)}d(x,y)^{\alpha
}\left( \varrho (x)^{\alpha }+\varrho (y)^{\alpha }\right) d\mu
(y)\right)
^{1/\alpha } \\
& \leq t\left( \varrho (x)+T_{t}^{\alpha }\varrho (x)\right) ,
\end{align*}
hence
\begin{eqnarray*}
E_{X^{(\alpha )}}(h,t) &\leq &t\left\| \varrho \right\| _{X^{(\alpha
)}}+t\left\| T_{t}^{\alpha }\varrho \text{\ }\right\| _{X^{(\alpha )}} \\
&\preceq &t\left\| \varrho \right\| _{X^{(\alpha )}}\text{ \ (by
Lemma \ref {Lfubi}).}
\end{eqnarray*}
In consequence
\begin{equation*}
E_{X^{(\alpha )}}(h,t)\preceq t\inf_{\varrho \in D(u)}\left\|
\varrho \right\| _{X^{(\alpha )}}=t\left\| h\right\|
_{\dot{M}^{1,X^{(\alpha )}}\left( \Omega \right) },
\end{equation*}
and taking the infimum over all representations of $f$ in $X^{(\alpha )}+%
\dot{M}^{1,X^{(\alpha )}}\left( \Omega \right) $, we obtain
\begin{equation*}
E_{X^{(\alpha )}}(f,t)\preceq K(f,t,X^{(\alpha
)},\dot{M}^{1,X^{(\alpha )}}).
\end{equation*}
For the converse inequality we will follow the same argument of the
proof of \cite[Theorem 4.1]{HIH}, with some elementary
modifications. For the reader's convenience, we will include here
the main steps. Let $f\in L^{\alpha }+L^{\infty \text{ }}$ and
$t>0$. By a standard covering argument, there is a covering of
$\Omega $ by balls $B_{i}=B(x_{i},t/6)$, $i\in \mathbb{N}$, such
that $\sum_{i}\chi _{2B_{i}}\leq N$ with the overlap constant $N>0$
depending only on $C_{\mu }$.

Let $\{\varphi _{i}\}_{i\in \mathbb{N}}$ be a collection of $Ct^{-1}-$%
Lipschitz functions $\varphi _{i}:\Omega \rightarrow \lbrack 0,1]$
such that
supp $\varphi _{i}\in $ $2B_{i}$ and $\sum_{i}\varphi _{i}(x)=1$ for all $%
x\in \Omega $. Let $h:\Omega \rightarrow \mathbb{R}$ defined by
\begin{equation*}
h(x)=\sum_{i\in N}m_{f}(B_{i})\varphi _{i}(x),\text{ \ \ }x\in
\Omega
\end{equation*}
where $m_{f}(B_{i})$ is the median value\footnote{%
The median value of a measurable function $u$ on a set $A\subset
\Omega $ is
\par
\begin{equation*}
m_{u}(A)=\max_{a\in \mathbb{R}}\left\{ \mu \left( \left\{ x\in
A:u(x)<a\right\} \right) \leq \frac{\mu (A)}{2}\right\} .
\end{equation*}
} of $f$ in $B_{i},$ and let $g=f-h,$ i.e.
\begin{equation*}
g(x)=\sum_{i\in N}\left( f(x)-m_{f}(B_{i})\right) \varphi
_{i}(x)=\sum_{i\in I_{x}}\left( f(x)-m_{f}(B_{i})\right) \varphi
_{i}(x)
\end{equation*}
where $I_{x}=\left\{ x:x\in 2B_{i}\right\} .$ The function $g$
satisfies (see \cite[formula (4.5)]{HIH})
\begin{equation}
\left| g(x)\right| \leq C\left( \frac{1}{\mu
(B(x,t))}\int_{B(x,t)}\left| f(x)-f(z)\right| ^{\alpha }d\mu
(z)\right) ^{1/\alpha }=C\nabla _{t}^{\alpha }f(x).  \label{ggg}
\end{equation}
Therefore
\begin{equation*}
\left\| g\right\| _{X^{(\alpha )}}\leq CE_{X^{(\alpha )}}(f,t).
\end{equation*}
Next we estimate $h$ in the $\dot{M}^{1,X^{(\alpha )}}$ norm. Let
$t>0$ and let $x,y\in \Omega $. We consider the following two cases.

Case 1: If $d(x,y)\leq t$, then (see \cite[page 352]{HIH})
\begin{equation*}
\left| h(x)-h(y)\right| \leq C\frac{d(x,y)}{t}\left( \frac{1}{\mu (B(x,2t))}%
\int_{B(x,2t)}\left| f(x)-f(z)\right| ^{\alpha }d\mu (z)\right)
^{1/\alpha }
\end{equation*}

Case 2: Let $d(x,y)>t.$ Since
\begin{equation*}
\left| h(x)-h(y)\right| \leq \left| f(x)-f(y)\right| +\left|
g(x)\right| +\left| g(y)\right| ,
\end{equation*}
it suffices to estimate the terms on the right side. The assumption $%
d(x,y)>t $ and (\ref{ggg}) imply that
\begin{equation*}
\left| g(x)\right| \leq C\frac{d(x,y)}{t}\left( \frac{1}{\mu (B(x,t))}%
\int_{B(x,t)}\left| f(x)-f(z)\right| ^{\alpha }d\mu (z)\right)
^{1/\alpha }
\end{equation*}
and a corresponding upper bound holds for $|g(y)|.$

By \cite[page 352]{HIH} and the doubling property we obtain
\begin{equation*}
\left| f(x)-f(y)\right| \leq Cd(x,y)\left(
f_{t}^{\#}(x)+f_{t}^{\#}(y)\right)
\end{equation*}
where
\begin{equation*}
f_{t}^{\#}(x)=\sup_{r\geq t}\frac{1}{r}\left( \frac{1}{\mu (B(x,r))}%
\int_{B(x,r)}\left| f(x)-f(z)\right| ^{\alpha }d\mu (z)\right)
^{1/\alpha }.
\end{equation*}
Collecting the above estimates, we get, in both cases, that
\begin{equation*}
\left| h(x)-h(y)\right| \leq Cd(x,y)\left(
f_{t}^{\#}(x)+f_{t}^{\#}(y)\right) .
\end{equation*}
Using the definition of $f_{t}^{\#}$ and the doubling property, we
have
\begin{equation*}
f_{t}^{\#}(x)\leq \frac{C}{t}\sum_{j=0}^{\infty }2^{-j}\left(
\frac{1}{\mu (B(x,2^{j}t))}\int_{B(x,2^{j}t)}\left| f(x)-f(z)\right|
^{\alpha }d\mu (z)\right) ^{1/\alpha }.
\end{equation*}
Since $0<\alpha \leq 1$, we obtain\footnote{\label{note1}recall that
$\left(
\sum_{i\in \mathbb{N}}a_{i}\right) ^{\alpha }\leq \sum_{i\in \mathbb{N}%
}a_{i}^{\alpha },$ if $a_{i}\geq 0$ and $0<\alpha \leq 1$.}
\begin{equation}
\left\| f_{t}^{\#}(x)\right\| _{X^{(\alpha )}}^{\alpha }\leq \frac{C}{%
t^{\alpha }}\sum_{j=0}^{\infty }2^{-j\alpha }\left\| \Delta
_{2^{j}t}^{\alpha }f(x)\right\| _{X^{(\alpha )}}^{\alpha }=\frac{C}{%
t^{\alpha }}\sum_{j=0}^{\infty }2^{-j\alpha }E_{X^{(\alpha
)}}(f,2^{j}t)^{\alpha }.  \label{ggg1}
\end{equation}
Thus, the required inequality
\begin{equation*}
K(f,t,X^{(\alpha )},\dot{M}^{1,X^{(\alpha )}})\leq C\left(
\sum_{k=0}^{\infty }2^{-j\alpha }E_{X^{(\alpha )}}(f,2^{j}t)^{\alpha
}\right) ^{1/\alpha }
\end{equation*}
is obtained using (\ref{ggg}), (\ref{ggg1}) and the definition of the $K-$%
functional.

Finally, equivalence (\ref{kf2}) follows with the same proof of
\cite[Theorem 4.1]{HIH}.
\end{proof}

\bigskip

Having determined the $K-$functional between $X^{(\alpha )}\ $and $\dot{M}%
^{1,X^{(\alpha )}}$ in terms of the $X^{(\alpha )}-$modulus of
smoothness, it is know routine to identify the corresponding
interpolation spaces.

\begin{corollary}
Suppose $X$ is a r.i. space and $0<\alpha <1.$ If $0<s<1$ and
$0<q\leq \infty ,$ then
\begin{equation*}
\dot{B}_{X^{(\alpha )},q}^{s}=\left( X^{(\alpha
)},\dot{M}^{1,X^{(\alpha )}}\right) _{s,q},
\end{equation*}
and
\begin{equation*}
B_{X^{(\alpha )},q}^{s}=\left( X^{(\alpha )},M^{1,X^{(\alpha
)}}\right) _{s,q},
\end{equation*}
with equivalent norms.
\end{corollary}

\begin{remark}
\label{equal}It follows from the previous result that if $1\leq
p<\infty ,$ then
\begin{equation*}
\dot{B}_{L^{p},q}^{s}=\left( L^{p},\dot{M}^{1,p}\right) _{s,q},
\end{equation*}
but also, by \cite[Theorem 4.1]{HIH}
\begin{equation*}
\mathcal{\dot{B}}_{p,q}^{s}=\left( L^{p},\dot{M}^{1,p}\right)
_{s,q},
\end{equation*}
thus
\begin{equation*}
\dot{B}_{L^{p},q}^{s}=\mathcal{\dot{B}}_{p,q}^{s}.
\end{equation*}
with equivalent norms.
\end{remark}

\section{The proofs of Theorems \ref{k1} and \ref{kk1}\label{section4}}

We will start with a lemma, which is a particular case of
\cite[Theorem 4] {MaOr1}, needed in the proof of our main results

\begin{lemma}
\label{teoMO1} Let $\left( \Omega ,d,\mu \right) $ be doubling
metric measure space with upper dimension $Q.$ Assume that there is
a constant $c>0$ such that
\begin{equation*}
\mu (B(x,r))\geq cr^{Q},\text{ \ }0<r\leq 1,\text{ }x\in \Omega .
\end{equation*}
Let $X$ be a r.i. space on $\Omega .$ Let $0<\alpha \leq 1,$ then for all $%
f\in \dot{M}^{1,X^{(\alpha )}}$ and $g\in D(f),$ we have that
\begin{equation*}
\left( f^{\alpha }\right) ^{\ast \ast }(t)-\left( f^{\alpha }\right)
^{\ast }(t)\leq ct^{\alpha /Q}\left( g^{\alpha }\right) ^{\ast \ast
}(t)\text{, \ \ \ }0<t<\frac{c}{2},
\end{equation*}
where $c=c(\alpha ,C_{\mu })$ is a constant that just depends on the
doubling constant $C_{\mu }$ and $\alpha $.
\end{lemma}

\begin{proof}
(\textbf{of Theorem \ref{k1}}) Given $X$ a r.i. space on $\Omega $ and
$0<\alpha \leq 1,$ we define
\begin{equation*}
K_{\alpha }(f,t,X^{(\alpha )},\dot{M}^{1,X^{(\alpha
)}})=\inf_{f=f_{0}+f_{1}}\left( \left\| f_{0}\right\| _{X^{(\alpha
)}}^{\alpha }+t^{\alpha }\left\| f_{1}\right\|
_{\dot{M}^{1,X^{(\alpha )}}}^{\alpha }\right) ^{1/\alpha }.
\end{equation*}
In what follows, for simplify, we shall write $K_{\alpha }(f,t,X^{(\alpha )},%
\dot{M}^{1,X^{(\alpha )}})=K_{\alpha }(f,t)$ and the classical $K-$%
functional $K(f,t,X^{(\alpha )},\dot{M}^{1,X^{(\alpha )}})$ will be
denoted by $K(f,t).$

(i)$\rightarrow$ (ii) Positivity will play a role in the arguments so it
will be useful to note for future use that if $\left\| \cdot
\right\| $ denotes either $\left\| \cdot \right\| _{X^{(\alpha )}}$
or $\left\| \cdot \right\| _{\dot{M}^{1,X^{(\alpha )}}},$ we have
\begin{equation}
\left\| \left| f\right| \right\| \leq \left\| f\right\| .
\label{dolores}
\end{equation}
Notice that if $h\in \dot{M}^{1,X^{(\alpha )}}$ and $\varrho \in
D(h),$ then $\left| h\right| \in \dot{M}^{1,X^{(\alpha )}}$ and
$\varrho \in D(\left| h\right| ),$ since
\begin{equation*}
\left| \left| h\right| (x)-\left| h\right| (y)\right| \leq \left|
h(x)-h(y)\right| \leq d(x,y)\left( \varrho (x)+\varrho (y)\right)
\end{equation*}
hence
\begin{equation*}
\left\| \left| h\right| \right\| _{\dot{M}^{1,X^{(\alpha
)}}}=\inf_{\varrho \in D(\left| h\right| )}\left\| \varrho \right\|
_{X^{(\alpha )}}\leq \inf_{\varrho \in D(h)}\left\| \varrho \right\|
_{X^{(\alpha )}}=\left\| h\right\| _{\dot{M}^{1,X^{(\alpha )}}}.
\end{equation*}
Let $\varepsilon >0$, and consider any decomposition $f=f-h+h,$ with
$h\in \dot{M}^{1,X^{(\alpha )}},$ such that
\begin{equation}
\left( \left\| f-h\right\| _{X^{(\alpha )}}^{\alpha }+t^{\alpha
}\left\| h\right\| _{\dot{M}^{1,X^{(\alpha )}}}^{\alpha }\right)
^{1/\alpha }\leq K_{\alpha }\left( f,t\right) +\varepsilon .
\label{k11}
\end{equation}
Since by (\ref{dolores}), $h\in \dot{M}^{1,X^{(\alpha )}}$ implies that $%
\left| h\right| \in \dot{M}^{1,X^{(\alpha )}}$, this decomposition
of $f$
produces the following decomposition of $\left| f\right| :$%
\begin{equation*}
\left| f\right| =\left| f\right| -\left| h\right| +\left| h\right| .
\end{equation*}
Therefore, by (\ref{dolores}) and (\ref{k11}) we have
\begin{eqnarray*}
\left( \left\| \left| f\right| -\left| h\right| \right\|
_{X^{(\alpha
)}}^{\alpha }+t^{\alpha }\left\| \left| h\right| \right\| _{\dot{M}%
^{1,X^{(\alpha )}}}^{\alpha }\right) ^{1/\alpha } &\leq &\left(
\left\|
f-h\right\| _{X^{(\alpha )}}^{\alpha }+t^{\alpha }\left\| h\right\| _{\dot{M}%
^{1,X^{(\alpha )}}}^{\alpha }\right) ^{1/\alpha } \\
&\leq &K_{\alpha }\left( f,t\right) +\varepsilon .
\end{eqnarray*}
Consequently
\begin{equation}
\inf_{0\leq h\in \dot{M}^{1,X^{(\alpha )}}}\left( \left\| \left|
f\right|
-h\right\| _{X^{(\alpha )}}^{\alpha }+t^{\alpha }\left\| h\right\| _{\dot{M}%
^{1,X^{(\alpha )}}}^{\alpha }\right) ^{1/\alpha }\leq K_{\alpha
}\left( f,t\right) .  \label{dolores1}
\end{equation}
Given $0\leq h\in \dot{M}^{1,X^{(\alpha )}}$ consider the
decomposition
\begin{equation*}
\left| f\right| =\left| f\right| -h+h.
\end{equation*}
Since $0<\alpha \leq 1,$ we get
\begin{equation*}
\left| f\right| ^{\alpha }=\left| \left| f\right| -h+h\right|
^{\alpha }\leq \left| \left| f\right| -h\right| ^{\alpha }+h^{\alpha
},
\end{equation*}
which by (\ref{readoble}) implies that
\begin{equation*}
\left( \left| f\right| ^{\alpha }\right) ^{\ast \ast }(t)\leq \left(
\left| \left| f\right| -h\right| ^{\alpha }\right) ^{\ast \ast
}(t)+\left( h^{\alpha }\right) ^{\ast \ast }(t).
\end{equation*}
On the other hand
\begin{equation*}
h=h-\left| f\right| +\left| f\right| \Rightarrow h^{\alpha }\leq
\left| f\right| ^{\alpha }+\left| \left| f\right| -h\right| ^{\alpha
},
\end{equation*}
and by (\ref{xxx}), we get
\begin{equation*}
\left( h^{\alpha }\right) ^{\ast }(2t)\leq \left( \left| f\right|
^{\alpha }+\left| \left| f\right| -h\right| ^{\alpha }\right) ^{\ast
}(2t)\leq f^{\ast }(t)^{\alpha }+\left( \left| \left| f\right|
-h\right| ^{\alpha }\right) ^{\ast }(t)\text{ }
\end{equation*}
hence
\begin{equation*}
f^{\ast }(t)^{\alpha }\geq \left( h^{\alpha }\right) ^{\ast
}(2t)-\left( \left| \left| f\right| -h\right| ^{\alpha }\right)
^{\ast }(t).
\end{equation*}
Combining the previous estimates we can write
\begin{align}
I(t)& :=\left( \left( \left| f\right| ^{\alpha }\right) ^{\ast \ast
}(t)-\left( \left| f\right| ^{\alpha }\right) ^{\ast }(t)\right)
\label{flaca} \\
& \leq \left( \left| \left| f\right| -h\right| ^{\alpha }\right)
^{\ast \ast }(t)+\left( h^{\alpha }\right) ^{\ast \ast }(t)-(\left(
h^{\alpha }\right) ^{\ast }(2t)-\left( \left| \left| f\right|
-h\right| ^{\alpha }\right)
^{\ast }(t))  \notag \\
& =\left( \left| \left| f\right| -h\right| ^{\alpha }\right) ^{\ast
\ast }(t)+\left( \left| \left| f\right| -h\right| ^{\alpha }\right)
^{\ast }(t))+\left( h^{\alpha }\right) ^{\ast \ast }(t)-\left(
h^{\alpha }\right)
^{\ast }(2t)  \notag \\
& \leq 2\left( \left| \left| f\right| -h\right| ^{\alpha }\right)
^{\ast \ast }(t)+\left( \left( h^{\alpha }\right) ^{\ast \ast
}(t)-\left( h^{\alpha
}\right) ^{\ast }(2t)\right)  \notag \\
& =2\left( \left| \left| f\right| -h\right| ^{\alpha }\right) ^{\ast
\ast }(t)+\left( \left( h^{\alpha }\right) ^{\ast \ast }(t)-\left(
h^{\alpha }\right) ^{\ast \ast }(2t)\right) +(\left( h^{\alpha
}\right) ^{\ast \ast
}(2t)-\left( h^{\alpha }\right) ^{\ast }(2t))  \notag \\
& =(I)+(II)+(III).  \notag
\end{align}
We first show that $(II)\leq (III)$ using fundamental theorem of
Calculus and (\ref{der2est}) to estimate $(II)$ as follows:
\begin{align*}
(II)& =\left( h^{\alpha }\right) ^{\ast \ast }(t)-\left( h^{\alpha
}\right)
^{\ast \ast }(2t) \\
& =\int_{t}^{2t}\left( \left( h^{\alpha }\right) ^{\ast \ast
}(s)-\left(
h^{\alpha }\right) ^{\ast }(s)\right) \frac{ds}{s} \\
& \leq 2t\left( \left( h^{\alpha }\right) ^{\ast \ast }(2t)-\left(
h^{\alpha
}\right) ^{\ast }(2t)\right) \int_{t}^{2t}\frac{ds}{s^{2}} \\
& =\left( \left( h^{\alpha }\right) ^{\ast \ast }(2t)-\left(
h^{\alpha
}\right) ^{\ast }(2t)\right) \\
& =(III).
\end{align*}
Inserting this information in (\ref{flaca}) we obtain
\begin{equation}
I(t)\leq 2\left( \left( \left| \left| f\right| -h\right| ^{\alpha
}\right) ^{\ast \ast }(t)+\left( h^{\alpha }\right) ^{\ast \ast
}(2t)-\left( h^{\alpha }\right) ^{\ast }(2t)\right) .
\label{flaca1}
\end{equation}

On the other hand, inequality (\ref{iteradobla}) with $R=1$ and
non-collapsing condition (\ref{noncolap}) imply the following lower
bound for the growth of the measure
\begin{equation*}
\mu (B(x,r))\geq \frac{b}{4^{Q}}r^{Q},\text{ \ }0<r\leq 1,
\end{equation*}
which by Lemma \ref{teoMO1} implies the existence of a positive constant $%
c=c(\alpha ,C_{\mu })\geq 1$ such that if $h\in
\dot{M}^{1,X^{(\alpha )}}$ and $g\in D(h),$ then
\begin{equation*}
\left( h^{\alpha }\right) ^{\ast \ast }(t)-\left( h^{\alpha }\right)
^{\ast }(t)\leq ct^{\alpha /Q}\left( g^{\alpha }\right) ^{\ast \ast
}(t)\text{, \ \ \ }0<t<\left( \frac{b}{2}\right) \frac{1}{4^{Q}}.
\end{equation*}
Inserting this information in (\ref{flaca1}) we have that if $0<t<\frac{b}{%
4^{Q+1}},$ then
\begin{align}
I(t)& \leq 2\left( \left( \left| \left| f\right| -h\right| ^{\alpha
}\right) ^{\ast \ast }(t)+ct^{\alpha /Q}\left( \left( g^{\alpha
}\right) ^{\ast \ast
}(t)\right) \right)  \label{det5} \\
& \leq 2c\left( \left( \left| \left| f\right| -h\right| ^{\alpha
}\right) ^{\ast \ast }(t)+t^{\alpha /Q}\left( \left( g^{\alpha
}\right) ^{\ast \ast
}(t)\right) \right)  \notag \\
& =2c\left( A(t)+B(t)\right) \text{.}  \notag
\end{align}
We now estimate the two terms on the right hand side of
(\ref{det5}). For
the term $A(t)$: Note that for any $G\in X^{(\alpha )},$%
\begin{equation*}
\left( \left| G\right| ^{\alpha }\right) ^{\ast \ast }(t)=\frac{1}{t}%
\int_{0}^{t}\left( \left| G\right| ^{\alpha }\right) ^{\ast }(s)ds=\frac{1}{t%
}\int_{0}^{1}\left( \left| G\right| ^{\alpha }\right) ^{\ast
}(s)\chi _{(0,t)}(s)ds.
\end{equation*}
Therefore, by H\"{o}lder's inequality (\ref{holin}) and (\ref{dual})
we have
\begin{align}
\left( \left| \left| f\right| -h\right| ^{\alpha }\right) ^{\ast \ast }(t)& =%
\frac{1}{t}\int_{0}^{1}\left( \left| \left| f\right| -h\right|
^{\alpha
}\right) ^{\ast }(s)\chi _{(0,t)}(s)ds  \label{det1} \\
& \leq \left\| \left( \left| f\right| -h\right) ^{\alpha }\right\| _{X}\frac{%
\phi _{X^{\prime }}(t)}{t}  \notag \\
& =\left\| \left( \left| f\right| -h\right) \right\| _{X^{(\alpha
)}}^{\alpha }\frac{1}{\phi _{X}(t)}.  \notag
\end{align}
Similarly, for $B(t)$ we get
\begin{equation}
B(t)=t^{\alpha /Q}\left( \left( g^{\alpha }\right) ^{\ast \ast
}(t)\right) \leq t^{\alpha /Q}\frac{\left\| g\right\| _{X^{(\alpha
)}}^{\alpha }}{\phi _{X}(t)}.  \label{det2}
\end{equation}
Inserting (\ref{det1}) and (\ref{det2}) back in (\ref{det5}) we find that$,$%
\begin{equation*}
I(t)\leq \frac{1}{\phi _{X}(t)}2c\left( \left\| \left( \left|
f\right| -h\right) \right\| _{X^{(\alpha )}}^{\alpha }+t^{\alpha
/Q}\left\| g\right\| _{X^{(\alpha )}}^{\alpha }\right) .
\end{equation*}
Thus
\begin{equation*}
I(t)^{1/\alpha }\leq \frac{\left( 2c\right) ^{\frac{1}{\alpha
}}}{\phi _{X^{(\alpha )}}(t)}\left( \left\| \left( \left| f\right|
-h\right) \right\|
_{X^{(\alpha )}}^{\alpha }+t^{\alpha /Q}\left\| h\right\| _{\dot{M}%
^{1,X^{(\alpha )}}}^{\alpha }\right) ^{1/\alpha }.
\end{equation*}
Therefore, by (\ref{dolores1})
\begin{eqnarray*}
I(t)^{1/\alpha } &\leq &\frac{\left( 2c\right) ^{\frac{1}{\alpha
}}}{\phi _{X^{(\alpha )}}(t)}\inf_{0\leq h\in \dot{M}^{1,X^{(\alpha
)}}}\left( \left\| \left( \left| f\right| -h\right) \right\|
_{X^{(\alpha )}}^{\alpha }+t^{\alpha /Q}\left\| h\right\|
_{\dot{M}^{1,X^{(\alpha )}}}^{\alpha
}\right) ^{1/\alpha } \\
&=&\left( 2c\right) ^{\frac{1}{\alpha }}\frac{K_{\alpha
}(f,t^{1/Q})}{\phi
_{X^{(\alpha )}}(t)} \\
&\leq &\left( 2c\right) ^{\frac{1}{\alpha }}2^{\frac{1}{\alpha }-1}\frac{%
K(f,t^{1/Q})}{\phi _{X^{(\alpha )}}(t)}\text{.}
\end{eqnarray*}
In summary we have proved that
\begin{equation}
O(\left| f\right| ^{\alpha },t)^{1/\alpha }\leq C\frac{K(f,t^{1/Q},X^{(%
\alpha )},\dot{M}^{1,X^{(\alpha )}})}{\phi _{X^{(\alpha )}}(t)},\text{ }0<t<%
\frac{b}{4^{Q+1}}.  \label{fundam4}
\end{equation}
Let $T=\min (\frac{b}{4^{Q+1}},1).$ By (\ref{fundam4}) we have
\begin{eqnarray*}
I &=&\int_{0}^{1}\left( O(\left| f\right| ^{\alpha },t)^{1/\alpha }\frac{%
\phi _{X^{(\alpha )}}(t)}{t^{s/Q}}\right) ^{q}\frac{dt}{t} \\
&=&\int_{0}^{T}\left( O(\left| f\right| ^{\alpha },t)^{1/\alpha
}\frac{\phi _{X^{(\alpha )}}(t)}{t^{s/Q}}\right)
^{q}\frac{dt}{t}+\int_{T}^{1}\left(
O(\left| f\right| ^{\alpha },t)^{1/\alpha }\frac{\phi _{X^{(\alpha )}}(t)}{%
t^{s/Q}}\right) ^{q}\frac{dt}{t} \\
&\leq &C\left( \int_{0}^{T}\left( \frac{K(f,t^{1/Q})}{t^{s/Q}}\right) ^{q}%
\frac{dt}{t}\right) +\int_{T}^{1}\left( \left( \left( \left|
f\right| ^{\alpha }\right) ^{\ast \ast }(t)\right) ^{1/\alpha
}(t)\frac{\phi
_{X^{(\alpha )}}(t)}{t^{s/Q}}\right) ^{q}\frac{dt}{t} \\
&\leq &C\left( \int_{0}^{\infty }\left(
\frac{K(f,t^{1/Q})}{t^{s/Q}}\right) ^{q}\frac{dt}{t}\right) +\left(
\left( \left| f\right| ^{\alpha }\right) ^{\ast \ast }(T)\right)
^{q/\alpha }\int_{T}^{1}\left( \frac{\phi
_{X^{(\alpha )}}(t)}{t^{s/Q}}\right) ^{q}\frac{dt}{t} \\
&\preceq &\int_{0}^{\infty }\left( \frac{K(f,z^{1/Q})}{z^{s/Q}}\right) ^{q}%
\frac{dz}{z}+\left\| f\right\| _{L^{\alpha }+L^{\infty }}^{q} \\
&=&\left\| f\right\| _{B_{X^{(\alpha )},q}^{s}}^{q}+\left\|
f\right\| _{L^{\min \alpha }+L^{\infty }}^{q},
\end{eqnarray*}
as we wanted to see.

(ii)$\rightarrow$(i) Let $0<s<\min (Q,1).$ Taking into account
(\ref{der2est}), by the fundamental Theorem of Calculus, we get
\begin{equation*}
f^{\ast \ast }(t)=\int_{t}^{1}\left( f^{\ast \ast }(z)-f^{\ast
}(z)\right) \frac{dz}{z}+f^{\ast \ast }(1),
\end{equation*}
consequently, by Fubini's theorem we get
\begin{align}
\int_{0}^{1}\frac{f^{\ast \ast }(t)}{t^{s/Q}}dt& \leq
\int_{0}^{1}\left( \int_{t}^{1}\left( f^{\ast \ast }(z)-f^{\ast
}(z)\right) \frac{dz}{z}\right) \frac{dt}{t^{s/Q}}+\frac{f^{\ast
\ast }(1)}{1-s/Q}\text{ (since }s/Q<1)
\notag \\
& =\int_{0}^{1}\left( f^{\ast \ast }(t)-f^{\ast }(t)\right) \left( \frac{1}{t%
}\int_{0}^{t}\frac{dz}{z^{s/Q}}\right) dt+\frac{f^{\ast \ast
}(1)}{1-s/Q)}
\notag \\
& =\frac{1}{1-s/Q}\left( \int_{0}^{1}\left( f^{\ast \ast
}(t)-f^{\ast
}(t)\right) \frac{dt}{t^{s/Q}}+f^{\ast \ast }(1)\right)  \notag \\
& =\frac{1}{1-s/Q}\left( \int_{0}^{1}\left( f^{\ast \ast
}(t)-f^{\ast }(t)\right) \frac{dt}{t^{s/Q}}+\left\| f\right\|
_{L^{1}+L^{\infty }}\right) .  \label{ups1}
\end{align}
Considering $X=L^{1},$ $\alpha =1,$ in (\ref{ineoscil}) we obtain
\begin{equation*}
f^{\ast \ast }(t)-f^{\ast }(t)\leq C\frac{K(f,t^{1/Q},L^{1},\dot{M}^{1,1})}{t%
},\text{ \ }0<t<T=\min (\frac{b}{4^{Q+1}},1)\text{ \ \ (since }\phi
_{L^{1}}(t)=t\text{).}
\end{equation*}
Inserting this information into (\ref{ups1}) we find that
\begin{align}
\int_{0}^{1}\frac{f^{\ast \ast }(t)}{t^{s/Q}}& \preceq \int_{0}^{T}\frac{%
K(f,t^{1/Q},L^{1},\dot{M}^{1,1})}{t^{s/Q}}\frac{dt}{t}+\int_{T}^{1}\left(
f^{\ast \ast }(t)-f^{\ast }(t)\right) \frac{dt}{t^{s/Q}}+\left\|
f\right\|
_{L^{1}+L^{\infty }}  \notag \\
& \preceq \int_{0}^{\infty }\frac{K(f,t^{1/Q},L^{1},\dot{M}^{1,1})}{z^{s}}%
\frac{dz}{z}+\int_{T}^{1}f^{\ast \ast }(t)\frac{dt}{t^{s/Q}}+\left\|
f\right\| _{L^{1}+L^{\infty }}  \notag \\
& \leq \int_{0}^{\infty }\frac{K(f,t^{1/Q},L^{1},\dot{M}^{1,1})}{z^{s}}\frac{%
dz}{z}+\int_{T}^{1}f^{\ast \ast }(t)\frac{dt}{t^{s/Q}}+\left\|
f\right\|
_{L^{1}+L^{\infty }}  \notag \\
& \preceq \left\| f\right\| _{\dot{B}_{1,1}^{s}\left( \Omega \right)
}+f^{\ast \ast }(T)\int_{T}^{1}\frac{dt}{t^{s/Q}}+\left\| f\right\|
_{L^{1}+L^{\infty }}  \notag \\
& \preceq \left\| f\right\| _{\dot{B}_{1,1}^{s}\left( \Omega \right)
}+\left\| f\right\| _{L^{1}+L^{\infty }}.  \label{ups0}
\end{align}

For a fixed $x_{0}\in \Omega ,$ we define the Lipschitz function
\begin{equation*}
u_{x_{0}}(y):=\left\{
\begin{array}{ll}
2-d(x_{0},y) & \text{if }y\in B(x_{0},2)\backslash B(x_{0},1), \\
1 & \text{if }y\in B(x_{0},1), \\
0 & \text{if }y\in \Omega \backslash B(x_{0},2).
\end{array}
\right.
\end{equation*}
It is easily seen that $g_{x_{0}}(y)=\chi _{B(x_{0},2)}(y)\in
D(u_{x_{0}})$

By Fubini's theorem
\begin{align}
E_{L^{1}}(u_{x_{0}},t)& \leq \int_{\Omega }\left|
u_{x_{0}}(x)\right| d\mu (x)+\int_{\Omega }\frac{1}{\mu
(B(x,t))}\int_{B(x,t)}\left|
u_{x_{0}}(y)\right| d\mu (y)d\mu (x)  \notag \\
& \leq \left\| u_{x_{0}}\right\| _{L^{1}}+\int_{\Omega }\left|
u_{x_{0}}(y)\right| \left( \int_{B(y,t)}\frac{1}{\mu (B(x,t))}d\mu
(x)\right) d\mu (y)  \notag \\
& \preceq \left\| u_{x_{0}}\right\| _{L^{1}},  \label{hello}
\end{align}
the last estimate follows from the doubling property of $\mu $ and since $%
B(y,t)\subset B(x,2t)$ whenever $x\in B(y,t).$

Using that $g_{x_{0}}\in D(u_{x_{0}})$ and with a similar argument as in (%
\ref{hello}), we get
\begin{align}
E_{L^{1}}(u_{x_{0}},t)& =\int_{\Omega }\left( \frac{1}{\mu (B(x,t))}%
\int_{B(x,t)}\left| u_{x_{0}}(x)-u_{x_{0}}(y)\right| d\mu (y)\right)
d\mu (x)
\notag \\
& \leq \int_{\Omega }\left( \frac{1}{\mu
(B(x,t))}\int_{B(x,t)}d(x,y)\left|
g_{x_{0}}(x)+g_{x_{0}}(y)\right| d\mu (y)\right) d\mu (x)  \notag \\
& \leq t\left( \int_{\Omega }\left| g_{x_{0}}(x)\right| d\mu
(x)+\int_{\Omega }\frac{1}{\mu (B(x,t))}\int_{B(x,t)}\left|
g_{x_{0}}(y)\right| d\mu (y)d\mu (x)\right)  \notag \\
& \preceq t\left\| g_{x_{0}}\right\| _{L^{1}}.  \label{hello1}
\end{align}
Thus, combining (\ref{hello}) and (\ref{hello1}), we see that
\begin{align*}
E_{L^{1}}(u_{x_{0}},t)& \preceq \min (\left\| u_{x_{0}}\right\|
_{L^{1}},t\left\| g_{x_{0}}\right\| _{L^{1}}) \\
& \leq \min (\mu \left( B(x_{0},2)\right) ,t\mu \left( B(x_{0},2)\right) ) \\
& \preceq \min (1,t)\mu \left( B(x_{0},2)\right) .
\end{align*}
Therefore,
\begin{equation*}
\left\| u_{x_{0}}\right\| _{\dot{B}_{1,1}^{s}}\preceq \mu \left(
B(x_{0},2)\right) ,
\end{equation*}
and obviously
\begin{equation*}
\left\| u_{x_{0}}\right\| _{L^{1}+L^{\infty }}\leq \left\|
u_{x_{0}}\right\| _{L^{1}}\leq \mu \left( B(x_{0},2)\right) .
\end{equation*}
On the other hand,
\begin{equation*}
\int_{0}^{1}\frac{\left( u_{x_{0}}\right) ^{\ast \ast
}(t)}{t^{s/Q}}\geq \int_{0}^{1}\frac{\left( u_{x_{0}}\right) ^{\ast
}(t)}{t^{s/Q}}\geq
\int_{0}^{\min (1,\mu \left( B(x_{0},1)\right) }\frac{dt}{t^{s/Q}}=\frac{%
\left( \min (1,\mu \left( B(x_{0},1)\right) \right)
^{1-s/Q}}{1-s/Q}.
\end{equation*}
The previous computation and (\ref{ups0}) implies that
\begin{equation*}
\frac{\left( \min (1,\mu \left( B(x_{0},1)\right) \right) ^{1-s/Q}}{1-s/Q}%
\preceq 2\mu \left( B(x_{0},2)\right) \leq 2C_{\mu }\mu \left(
B(x_{0},1)\right) ,
\end{equation*}
and from here we conclude that
\begin{equation*}
1\preceq \mu \left( B(x_{0},1)\right)
\end{equation*}
which completes the proof.
\end{proof}

\begin{remark}
Theorem \ref{kk1} follows from Theorem \ref{k1} considering the r.i. spaces $%
X^{(1/\alpha )}.$
\end{remark}

\begin{remark}
Inequality (\ref{fundam4}) fits into the research program developed in Mart%
\'{i}n and Milman where Sobolev-Poincar\'{e} inequalities for
Lipschitz functions and $K-$functionals are widely considered (for
further details, we refer the reader to \cite{MM4}, \cite{Mar} and
\cite{Mas}).
\end{remark}

\section{Embedding Theorems for Haj\l asz-Besov spaces into q.r.i. spaces%
\label{sec05}}

In this section is devoted to find Sobolev type embedding for $\dot{B}%
_{X^{(\alpha )},q}^{s}\left( \Omega \right) $ spaces into a q.r.i.
spaces. The proofs make substantial use of the (nonlinear) spaces
$L(v,\alpha ,q)$
defined to be the set of all functions $f\in L^{\alpha }+L^{\infty }$ ($%
0<\alpha \leq 1$) such that
\begin{equation*}
\left\| f\right\| _{L(v,\alpha ,q)}=\left\{
\begin{array}{ll}
\left( \dint_{0}^{1}\left( O(\left| f\right| ^{\alpha },z)^{1/\alpha
}v(z)\right) ^{q}\frac{dz}{z}\right) ^{1/q}, & 0<q<\infty , \\
\sup_{0<t<1}\left( O(\left| f\right| ^{\alpha },t)\right) v(t), &
q=\infty .
\end{array}
\right.
\end{equation*}
is finite. (Here $v$ is a weight (a positive locally integrable function on $%
\left( 0,\infty \right) )$.

Before continuing, we note $\left\| f\right\| _{L(v,\alpha ,q)}$
depends neither on the growth of $(|f|^{\alpha })^{\ast }$ nor on
$(|f|^{\alpha })^{\ast \ast }$ but rather on the oscillation, which
causes obstacles for applications, to overcome such difficulties we
shall need the next results.

\begin{remark}
\label{oscil}Given $q$ a positive number and $f$ a $\mu -$measurable
function on $\Omega $, it follows ready that
\begin{equation*}
\frac{\partial }{\partial t}\left( \left| f\right| ^{\alpha }\right)
^{\ast \ast }(t)^{q}=-q\left( \left| f\right| ^{\alpha }\right)
^{\ast \ast }(t)^{q-1}O(\left| f\right| ^{\alpha },t)\frac{1}{t},
\end{equation*}
therefore, by the Fundamental Theorem of Calculus,
\begin{eqnarray}
\left( \left| f\right| ^{\alpha }\right) ^{\ast \ast }(t)^{q}
&=&q\left( \int_{t}^{1}\left( \left| f\right| ^{\alpha }\right)
^{\ast \ast }(t)^{q-1}\left( O(\left| f\right| ^{\alpha },t)\right)
\frac{dt}{t}+\left(
\left| f\right| ^{\alpha }\right) ^{\ast \ast }(1)^{q}\right) \   \notag \\
&=&q\int_{t}^{1}\left( \left| f\right| ^{\alpha }\right) ^{\ast \ast
}(t)^{q-1}\left( O(\left| f\right| ^{\alpha },t)\right) \frac{dt}{t}%
+q\left\| f\right\| _{L^{^{\alpha }}+L^{\infty }}^{q}.
\label{osci1}
\end{eqnarray}
In particular, if $0<q\leq 1,$ in view of $O(\left| f\right|
^{\alpha },t)\leq \left( \left| f\right| ^{\alpha }\right) ^{\ast
\ast }(t),$ we get
\begin{equation}
\left( \left| f\right| ^{\alpha }\right) ^{\ast \ast }(t)^{q}\leq
q\int_{t}^{1}\left( O(\left| f\right| ^{\alpha },z)\right) ^{q}\frac{dz}{z}%
+q\left\| f\right\| _{L^{^{\alpha }}+L^{\infty }}^{q}.
\label{osci2}
\end{equation}
\end{remark}

\begin{lemma}
\label{infinito}Let $v$ be a weight, $0<\alpha \leq 1\ $and $0<q\leq
\infty , $ then the following embedding holds
\begin{equation*}
\left\| f\right\| _{L^{\infty }}\preceq \left\| f\right\|
_{L(v,\alpha ,q)}+\left\| f\right\| _{L^{^{\alpha }}+L^{\infty }},
\end{equation*}
if $m_{v,\alpha ,q}(0)<\infty ,$ where
\begin{equation*}
m_{v,\alpha ,q}(t):=m(t)=\left\{
\begin{array}{ll}
\int_{t}^{1}\left( \frac{1}{v(z)}\right) ^{\frac{\alpha q}{q-\alpha }}\frac{%
dz}{z} & \text{if \ \ }\alpha <q\leq \infty , \\
\int_{t}^{1}\frac{1}{v(z)^{\alpha }}\frac{dz}{z} & \text{if \ \
}q=\infty ,
\\
\sup\limits_{z\in \lbrack t,1)}\frac{1}{v(z)} & \text{if \ \
}0<q\leq \alpha .
\end{array}
\right.
\end{equation*}
\end{lemma}

\begin{proof}
It follows from (\ref{osci1}) that
\begin{eqnarray}
\left\| f\right\| _{L^{\infty }}^{\alpha } &=&\left( \left| f\right|
^{\alpha }\right) ^{\ast \ast }(0)=\int_{0}^{1}\left( O(\left|
f\right| ^{\alpha },t)\right) \frac{dt}{t}+\left\| f\right\|
_{L^{^{\alpha
}}+L^{\infty }}^{\alpha }  \notag \\
&=&I+\left\| f\right\| _{L^{^{\alpha }}+L^{\infty }}^{\alpha }.
\notag
\end{eqnarray}
\textbf{Case }$\alpha <q<\infty .$ Since $\alpha /q>1$, by Holder
inequality we obtain
\begin{equation}
I\leq \left( \int_{0}^{1}\left( O(\left| f\right| ^{\alpha
},t)^{1/\alpha }v(t)\right) ^{q}\frac{dt}{t}\right) ^{\frac{\alpha
}{q}}\left(
\int_{0}^{1}\left( \frac{1}{v(t)}\right) ^{\frac{\alpha q}{q-\alpha }}\frac{%
dz}{t}\right) ^{\frac{q-\alpha }{q}},  \notag
\end{equation}
whence
\begin{equation*}
\left\| f\right\| _{L^{\infty }}\preceq \left( m(0)\right)
^{\frac{q-\alpha }{q\alpha }}\left\| f\right\| _{L(v,\alpha
,q)}+\left\| f\right\| _{L^{^{\alpha }}+L^{\infty }}.
\end{equation*}
\textbf{Case }$q=\infty .$ Obviously
\begin{equation*}
I\leq \sup_{t\in \lbrack 0,1]}\left( O(\left| f\right| ^{\alpha
},t)v^{\alpha }(t)\right) \int_{0}^{1}\frac{1}{v^{\alpha
}(t)}\frac{dt}{t},
\end{equation*}
thus
\begin{equation*}
\left\| f\right\| _{L^{\infty }}\preceq m(0)\left\| f\right\|
_{L(v,\alpha ,\infty )}+\left\| f\right\| _{L^{^{\alpha }}+L^{\infty
}}
\end{equation*}
\textbf{Case }$q\leq \alpha .$ By inequality (\ref{osci2}), we get
\begin{eqnarray*}
\left\| f\right\| _{L^{\infty }}^{q} &=&\left( \left| f\right|
^{\alpha }\right) ^{\ast \ast }(0)^{q}\leq \frac{q}{\alpha
}\int_{0}^{1}\left(
O(\left| f\right| ^{\alpha },t)^{1/\alpha }v(t)\right) ^{q}\frac{dt}{%
v(t)^{q}t}+\left\| f\right\| _{L^{^{\alpha }}+L^{\infty }}^{q} \\
&\leq &\frac{q}{\alpha }\sup_{t\in \lbrack 0,1]}\left( \frac{1}{v(z)^{q}}%
\right) \left\| f\right\| _{L(v,\alpha ,\infty )}^{q}+\left\|
f\right\| _{L^{^{\alpha }}+L^{\infty }}^{q},
\end{eqnarray*}
hence
\begin{equation*}
\left\| f\right\| _{L^{\infty }}\preceq m(0)\left\| f\right\|
_{L(v,\alpha ,q)}+\left\| f\right\| _{L^{^{\alpha }}+L^{\infty }},
\end{equation*}
as we wished to show.
\end{proof}

\begin{lemma}
\label{pesos}Let $v$ be a weight. Let $0<\alpha \leq 1$ and $0<q\leq
\infty $ be such that $m_{v,\alpha ,q}(0)=\infty ,$ then

\begin{enumerate}
\item  If $\alpha <q<\infty $, then
\begin{equation*}
\left( \int_{0}^{1}\left( \left( \left| f\right| ^{\alpha }\right)
^{\ast \ast }(t)^{1/\alpha }w(z)\right) ^{q}\frac{dz}{z}\right)
^{1/q}\preceq \left\| f\right\| _{L(v,\alpha ,q)}+\left\| f\right\|
_{L^{\alpha }+L^{\infty }},
\end{equation*}
where $w$ is defined by
\begin{equation*}
\frac{w^{q}(t)}{t}=\frac{\partial }{\partial t}\left(
1+\int_{t}^{1}\left( \frac{1}{v(z)}\right) ^{\frac{\alpha
q}{q-\alpha }}\frac{dz}{z}\right) ^{1-q/\alpha }.
\end{equation*}

\item  If $q=\infty ,$ then
\begin{equation*}
\sup_{t\in \lbrack 0,1]}\frac{\left( \left| f\right| ^{\alpha
}\right)
^{\ast \ast }(t)^{1/\alpha }}{\left( 1+\int_{t}^{1}\frac{1}{v^{\alpha }(z)}%
\frac{dz}{z}\right) ^{1/\alpha }}\preceq \left\| f\right\|
_{L(v,\alpha ,\infty )}
\end{equation*}

\item  If $0<q\leq \alpha ,$ then for any weight $u$ such that
\begin{equation*}
\int_{0}^{t}u^{q}(z)\frac{dz}{z}\preceq v^{q}(t),\text{ \ }0<t<1
\end{equation*}
we have that
\begin{equation*}
\left( \int_{0}^{1}\left( \left( \left| f\right| ^{\alpha }\right)
^{\ast \ast }(t)^{1/\alpha }u(t)\right) ^{q}\frac{dt}{t}\right)
^{1/q}\preceq \left\| f\right\| _{L(v,\alpha ,q)}+\left\| f\right\|
_{L^{\alpha }+L^{\infty }}.
\end{equation*}
\end{enumerate}
\end{lemma}

\begin{proof}
(i) By (\ref{osci1}) we have that
\begin{align}
I& =\int_{0}^{1}\left( \left| f\right| ^{\alpha }\right) ^{\ast \ast
}(t)^{q/\alpha }w^{q}(t)\frac{dt}{t}  \notag \\
& =\frac{q}{\alpha }\left( \int_{0}^{1}\int_{t}^{1}\left( \left|
f\right| ^{\alpha }\right) ^{\ast \ast }(z)^{q/\alpha -1}\left(
O(\left| f\right|
^{\alpha },z)\right) \frac{dz}{z}\right) w^{q}(t)\frac{dt}{t}  \notag \\
& \text{ \ \ \ \ }+\left( \left| f\right| ^{\alpha }\right) ^{\ast
\ast
}(1)^{q/\alpha }\int_{0}^{1}w^{q}(t)\frac{dt}{t}  \notag \\
& =\frac{q}{\alpha }II+\left( \left| f\right| ^{\alpha }\right)
^{\ast \ast }(1)^{q/\alpha }\int_{0}^{1}w^{q}(t)\frac{dt}{t}.
\label{222}
\end{align}
Applying Fubini's Theorem and H\"{o}lder inequality, we obtain
\begin{align}
II& =\int_{0}^{1}\left( \left| f\right| ^{\alpha }\right) ^{\ast \ast }(t)^{%
\frac{q}{\alpha }-1}\left( O(\left| f\right| ^{\alpha },t)\right)
\left(
\frac{1}{t}\int_{0}^{t}w^{q}(z)\frac{dz}{z}\right) dt  \notag \\
& =\int_{0}^{1}\left( \left| f\right| ^{\alpha }\right) ^{\ast \ast }(t)^{%
\frac{q}{\alpha }-1}\left( \frac{w^{q}(t)}{t}\right) ^{1-\frac{\alpha }{q}%
}\left( O(\left| f\right| ^{\alpha },t)\right) \left( \frac{w^{q}(t)}{t}%
\right) ^{\frac{\alpha }{q}-1}\left( \frac{1}{t}\int_{0}^{t}w^{q}(z)\frac{dz%
}{z}\right) dt  \notag \\
& \leq \left( \int_{0}^{1}\left( \left( \left| f\right| ^{\alpha
}\right) ^{\ast \ast }(t)w^{q}(t)\right) \frac{dt}{t}\right)
^{1-\frac{\alpha }{q}}
\notag \\
& \text{ \ \ \ \ }\times \left( \int_{0}^{1}\left( O(\left| f\right|
^{\alpha },t)\right) ^{\frac{q}{\alpha }}\left(
\frac{w^{q}(t)}{t}\right)
^{1-\frac{q}{\alpha }}\left( \frac{1}{t}\int_{0}^{t}w^{q}(z)\frac{dz}{z}%
\right) ^{\frac{q}{\alpha }}dt\right) ^{\frac{\alpha }{q}}
\label{222a}
\end{align}
On the other hand, an elementary computation shows that
\begin{equation*}
\left( \frac{w^{q}(t)}{t}\right) ^{1-\frac{q}{\alpha }}\left( \frac{1}{t}%
\int_{0}^{t}w^{q}(z)\frac{dz}{z}\right) ^{\frac{q}{\alpha }}\simeq \frac{%
v^{\alpha }(t)}{t},
\end{equation*}
whence combining (\ref{222}) and (\ref{222a}) we obtain
\begin{equation*}
\left( \int_{0}^{1}\left( \left| f\right| ^{\alpha }\right) ^{\ast
\ast
}(t)^{\frac{q}{\alpha }}w^{q}(t)\frac{dt}{t}\right) ^{\frac{\alpha }{q}%
}\preceq \left\| f\right\| _{L(v,\alpha ,q)}^{\alpha }+\frac{\left(
\left|
f\right| ^{\alpha }\right) ^{\ast \ast }(1)^{\frac{q}{\alpha }%
}\int_{0}^{1}w^{q}(t)\frac{dt}{t}}{\left( \int_{0}^{1}\left( \left|
f\right|
^{\alpha }\right) ^{\ast \ast }(t)^{\frac{q}{\alpha }}w^{q}(t)\frac{dt}{t}%
\right) ^{1-\frac{\alpha }{q}}}.
\end{equation*}
Finally, since
\begin{eqnarray*}
\frac{\left( \left| f\right| ^{\alpha }\right) ^{\ast \ast }(1)^{\frac{q}{%
\alpha }}\int_{0}^{1}w^{q}(t)\frac{dt}{t}}{\left( \int_{0}^{1}\left(
\left|
f\right| ^{\alpha }\right) ^{\ast \ast }(t)^{\frac{q}{\alpha }}w^{q}(t)\frac{%
dt}{t}\right) ^{1-\frac{\alpha }{q}}} &\leq &\frac{\left( \left|
f\right|
^{\alpha }\right) ^{\ast \ast }(1)^{\frac{q}{\alpha }}\int_{0}^{1}w^{q}(t)%
\frac{dt}{t}}{\left( \left( \left| f\right| ^{\alpha }\right) ^{\ast
\ast
}(1)^{\frac{q}{\alpha }}\int_{0}^{1}{}w^{q}(t)\frac{dt}{t}\right) ^{1-\frac{%
\alpha }{q}}} \\
&=&\left( \left| f\right| ^{\alpha }\right) ^{\ast \ast }(1)\left(
\int_{0}^{1}w^{q}(t)\frac{dt}{t}\right) ^{\frac{\alpha }{q}} \\
&=&\left\| f\right\| _{L^{^{\alpha }}+L^{\infty }}^{\alpha }\left(
\int_{0}^{1}w^{q}(t)\frac{dt}{t}\right) ^{\frac{\alpha }{q}}
\end{eqnarray*}
the desired inequality follows.

(ii) Using (\ref{osci1}) we can write
\begin{eqnarray*}
\left( \left| f\right| ^{\alpha }\right) ^{\ast \ast }(t)
&=&\int_{t}^{1}\left( O(\left| f\right| ^{\alpha },z)\right) \frac{dz}{z}%
+\left\| f\right\| _{L^{^{\alpha }}+L^{\infty }}^{\alpha } \\
&\leq &\sup_{t\in \lbrack 0,1]}\left( O(\left| f\right| ^{\alpha
},t)v^{\alpha }(t)\right) \int_{t}^{1}\frac{1}{v^{\alpha }(z)}\frac{dz}{z}%
+\left\| f\right\| _{L^{^{\alpha }}+L^{\infty }}^{\alpha } \\
&\leq &\left( \left\| f\right\| _{L(v,\alpha ,\infty )}^{\alpha
}+\left\| f\right\| _{L^{^{\alpha }}+L^{\infty }}^{\alpha }\right)
\left( 1+\int_{t}^{1}\frac{1}{v^{\alpha }(z)}\frac{dz}{z}\right) .
\end{eqnarray*}
Therefore
\begin{equation*}
\sup_{t\in \lbrack 0,1]}\frac{\left( \left| f\right| ^{\alpha
}\right)
^{\ast \ast }(t)^{1/\alpha }}{\left( 1+\int_{t}^{1}\frac{1}{v^{\alpha }(z)}%
\frac{dz}{z}\right) ^{1/\alpha }}\preceq \left\| f\right\|
_{L(v,\alpha ,\infty ).}
\end{equation*}

(iii) By (\ref{osci2}) and Fubini's Theorem we get
\begin{align*}
I& =\int_{0}^{1}\left( \left( \left| f\right| ^{\alpha }\right)
^{\ast \ast
}(t)^{1/\alpha }u(t)\right) ^{q}\frac{dt}{t} \\
& \leq \frac{q}{\alpha }\int_{0}^{1}\left( \int_{t}^{1}\left( O_{\mu
}(\left| f\right| ^{\alpha },z)\right) ^{q/\alpha
}\frac{dz}{z}\right) u^{q}(t)\frac{dt}{t}+\left\| f\right\|
_{L^{\alpha }+L^{\infty
}}^{q}\int_{0}^{1}u^{q}(t)\frac{dt}{t} \\
& =\frac{q}{\alpha }\int_{0}^{1}\left( O_{\mu }(\left| f\right|
^{\alpha
},t)\right) ^{q/\alpha }\left( \frac{1}{t}\int_{0}^{t}u^{q}(z)\frac{dt}{z}%
\right) dt+\left\| f\right\| _{L^{\alpha }+L^{\infty
}}^{q}\int_{0}^{1}u^{q}(t)\frac{dt}{t} \\
& \leq c\frac{q}{\alpha }\left\| f\right\| _{L(v,\alpha ,\infty
)}^{q}+\left\| f\right\| _{L^{\alpha }+L^{\infty }}^{q}\int_{0}^{1}u^{q}(t)%
\frac{dt}{t},
\end{align*}
which completes the proof.
\end{proof}

Now we are ready to proof the following embedding result.

\begin{theorem}
\label{embteo}Let $\left( \Omega ,d,\mu \right) $ be doubling with
upper dimension $Q$ that satisfies the non-collapsing condition. Let
$X$ be a r.i. space on $\Omega $, $0<\alpha \leq 1,$ $0<s<1,$
$0<q\leq \infty \ $and consider the function
\begin{equation*}
m_{X,s,\alpha ,q}:[0,1]\rightarrow \lbrack 0,\infty )
\end{equation*}
defined by
\begin{equation*}
m_{X,s,\alpha ,q}(t):=m(t)=\left\{
\begin{array}{ll}
\int_{t}^{1}\left( \frac{t^{\frac{s}{Q}}}{\phi _{X^{(\alpha )}}(z)}\right) ^{%
\frac{\alpha q}{q-\alpha }}\frac{dz}{z} & \text{if \ \ }\alpha
<q\leq \infty
, \\
\int_{t}^{1}\left( \frac{t^{\frac{s}{Q}}}{\phi _{X^{(\alpha
)}}}\right)
^{\alpha }\frac{dz}{z} & \text{if \ \ }q=\infty , \\
\sup\limits_{z\in \lbrack t,1)}\frac{t^{\frac{s}{Q}}}{\phi
_{X^{(\alpha )}}} & \text{if \ \ }0<q\leq \alpha .
\end{array}
\right.
\end{equation*}
Then

\begin{enumerate}
\item  If $m(0)<\infty ,$ then
\begin{equation*}
\left\| f\right\| _{L^{\infty }}\preceq \left\| f\right\|
_{B_{X^{(\alpha )},q}^{s}}+\left\| f\right\| _{L^{^{\alpha
}}+L^{\infty }}.
\end{equation*}

\item  If $m(0)=\infty $, then

\begin{enumerate}
\item  If $\alpha <q<\infty ,$ then
\begin{equation*}
\left( \int_{0}^{1}f^{\ast }(t)^{q}\frac{\left| m^{\prime }(t)\right| }{%
\left( 1+m(t)\right) ^{q/\alpha }}ds\right) ^{1/q}\preceq \left\|
f\right\| _{\dot{B}_{X^{(\alpha ),q}}^{s}}+\left\| f\right\|
_{L^{\alpha }+L^{\infty }.}
\end{equation*}

\item  If $q=\infty ,$ then
\begin{equation*}
\sup_{t\in \lbrack 0,1]}\left( \frac{f^{\ast }(t)^{q}}{\left(
1+m(t)\right) ^{1/\alpha }}\right) \preceq \left\| f\right\|
_{\dot{B}_{X^{(\alpha ),\infty }}^{s}}+\left\| f\right\| _{L^{\alpha
}+L^{\infty }.}
\end{equation*}

\item  If $0<q\leq \alpha ,$ then for any weight $u$ such that
\begin{equation*}
\int_{0}^{t}u^{q}(z)dz\preceq v^{q}(t),\text{ \ }0<t<1
\end{equation*}
there is a positive constant $C$ such that
\begin{equation*}
\left( \int_{0}^{1}\left( f^{\ast }(t)u(t)\right)
^{q}\frac{dt}{t}\right) ^{1/q}\leq C\left( \left\| f\right\|
_{\dot{B}_{X^{(\alpha ),q}}^{s}}+\left\| f\right\| _{L^{\alpha
}+L^{\infty }}\right) .
\end{equation*}
\end{enumerate}
\end{enumerate}
\end{theorem}

\begin{proof}
Let us write $K(f,t^{1/Q})=K(f,t^{1/Q},X^{(\alpha
)},\dot{M}^{1,X^{(\alpha )}})$ and $T=$ $\min (b/4^{Q+1},1)$. If
$T\leq 1,$ then
\begin{eqnarray*}
I &=&\left( \int_{0}^{T}\left( O(\left| f\right| ^{\alpha
},t)^{1/\alpha
}\left( \frac{\phi _{X}(t)^{1/\alpha }}{t^{\frac{s}{Q}}}\right) \right) ^{q}%
\frac{dt}{t}\right) ^{1/q} \\
&\leq &\left( \int_{0}^{1}\left( O(\left| f\right| ^{\alpha
},t)^{1/\alpha
}\left( \frac{\phi _{X^{(\alpha )}}(t)}{t^{\frac{s}{Q}}}\right) \right) ^{q}%
\frac{dt}{t}\right) ^{1/q}\text{\ (since }\phi _{X^{(\alpha
)}}=\left( \phi
_{X}\right) ^{1/\alpha }\text{)} \\
&\preceq &\left( \int_{0}^{\infty }\left( \frac{K(f,t^{1/Q})}{t^{s/Q}}%
\right) ^{q}\frac{dt}{t}\right) ^{1/q}\text{\ \ (by Theorem \ref{k1})} \\
&\preceq &\left\| f\right\| _{\dot{B}_{X^{(\alpha ),q}}^{s}}+\left\|
f\right\| _{L^{\alpha }+L^{\infty }}.
\end{eqnarray*}
If $T>1,$ then
\begin{eqnarray*}
I &=&\left( \int_{0}^{1}\left( O(\left| f\right| ^{\alpha
},t)^{1/\alpha
}\left( \frac{\phi _{X^{(\alpha )}}}{t^{\frac{s}{Q}}}(t)\right) \right) ^{q}%
\frac{dt}{t}\right) ^{1/q} \\
&\preceq &\left( \int_{0}^{T}\left( O(\left| f\right| ^{\alpha
},t)^{1/\alpha }\left( \frac{\phi _{X}(t)^{1/\alpha }}{t^{\frac{s}{Q}}}%
\right) \right) ^{q}\frac{dt}{t}\right) ^{1/q} \\
&&+\left( \int_{T}^{1}\left( O(\left| f\right| ^{\alpha
},t)^{1/\alpha
}\left( \frac{\phi _{X^{(\alpha )}}(t)}{t^{\frac{s}{Q}}}\right) \right) ^{q}%
\frac{dt}{t}\right) ^{1/q} \\
&=&(A)+(B)\text{.}
\end{eqnarray*}
The term $(A)$ is controlled as above, and
\begin{eqnarray*}
(B) &=&\left( \int_{T}^{1}\left( O(\left| f\right| ^{\alpha
},t)^{1/\alpha
}\left( \frac{\phi _{X^{(\alpha )}}(t)}{t^{\frac{s}{Q}}}\right) \right) ^{q}%
\frac{dt}{t}\right) ^{1/q} \\
&\leq &\left( \left| f\right| ^{\alpha }\right) ^{\ast \ast
}(T)^{1/\alpha
}\left( \int_{T}^{1}\left( \frac{\phi _{X^{(\alpha )}}(t)}{t^{\frac{s}{Q}}}%
\right) ^{q}\frac{dt}{t}\right) ^{1/q} \\
&\preceq &\left\| f\right\| _{L^{\alpha }+L^{\infty }}.
\end{eqnarray*}
In summary, in both cases we have proved that
\begin{equation*}
\left( \int_{0}^{1}\left( O(\left| f\right| ^{\alpha },t)^{1/\alpha
}\left(
\frac{\phi _{X^{(\alpha )}}}{t^{\frac{s}{Q}}}\right) \right) ^{q}\frac{dt}{t}%
\right) ^{1/q}\preceq \left\| f\right\| _{\dot{B}_{X^{(\alpha
),q}}^{s}}+\left\| f\right\| _{L^{\alpha }+L^{\infty }},
\end{equation*}
therefore, the result follows applying Lemmas \ref{infinito} and
\ref{pesos} to the weight
\begin{equation*}
v(t)=\frac{t^{\frac{s}{Q}}}{\phi _{X^{(\alpha )}}(t)}
\end{equation*}
and using that $f^{\ast }(t)\leq \left( \left| f\right| ^{\alpha
}\right) ^{\ast \ast }(t)^{1/\alpha }.$
\end{proof}

\begin{remark}
It is plain that if we work with Haj\l asz-Besov spaces
$B_{X^{(\alpha
)},q}^{s}$ instead of homogeneous Haj\l asz-Besov spaces $\dot{B}%
_{X^{(\alpha )},q}^{s}$, then Theorem \ref{embteo} remains true,
considering in the right hand side of the inequalities the term
$\left\| f\right\|
_{B_{X^{(\alpha ),q}}^{s}}$ instead of $\left\| f\right\| _{\dot{B}%
_{X^{(\alpha ),q}}^{s}}+\left\| f\right\| _{L^{\alpha }+L^{\infty
}}$, since by (\ref{embbb}) we know that $\left\| f\right\|
_{L^{\alpha }+L^{\infty }}\preceq \left\| f\right\| _{X^{(\alpha
)}}.$
\end{remark}

\subsubsection{Example: Homogeneous Lorentz-Zygmund-Haj\l asz-Besov spaces}

Lorentz-Zygmund spaces were introduced in \cite{BeRu}, they contain
many interesting nontrivial function spaces as Lebesgue spaces,
Lorentz spaces or Zygmund classes, which have important
applications, mainly in various limiting or critical situations, see
for example \cite{BreWa}. In what follows we investigate Haj\l
asz-Besov embedding results for such spaces.

Given $0<p<\infty ,$ $0<r\leq \infty $ and $\beta \in \mathbb{R},$
the
Lorentz-Zygmund space $L^{p,r}(\log L)^{\beta }$ consists of all $\mu -$%
measurable functions $f$ on $\Omega $ for which the quasi-norm
\begin{equation*}
\left\| f\right\| _{L^{p,r}(\log L)^{\beta }}=\left\{
\begin{array}{ll}
\left( \dint_{0}^{\infty }\left[ t^{\frac{1}{p}}\left( 1+\ln ^{+}\frac{1}{t}%
\right) ^{\beta }f^{\ast }(t)\right] ^{r}\frac{dt}{t}\right) ^{1/r},
&
0<r<\infty , \\
\sup_{t>0}t^{\frac{1}{p}}(1+\ln ^{+}\frac{1}{t})^{\beta }f^{\ast
}(t), & r=\infty ,
\end{array}
\right.
\end{equation*}
is finite. The fundamental function of $L^{p,r}(\log L)^{\beta }$
satisfies that
\begin{equation}
\phi _{L^{p,r}(\log L)^{\beta }}(t)\simeq t^{\frac{1}{p}}(1+\ln \frac{1}{t}%
)^{\beta }.  \label{funda}
\end{equation}

Lorentz-Zygmund spaces are $\alpha -$convex, indent taking into
account that if $1<p\leq \infty $ and $1\leq r\leq \infty ,$ the
functional $\left\|
\cdot \right\| _{L^{p,r}(\log L)^{\beta }}$ is equivalent to a norm, and $%
L^{1,1}(\log L)^{\beta }$ is a Banach space if $\beta \geq 0,$ it
follows ready that choosing $0<\alpha \leq 1$ satisfying
\begin{equation}
\left\{
\begin{array}{ll}
\alpha <p & \text{if }0<p<r\text{ with }0<p\leq 1,\text{ or }0<p=r\leq 1%
\text{ and }\beta <0, \\
\alpha =\min (1,r) & \text{otherwise,}
\end{array}
\right.  \label{elealpha}
\end{equation}
we have that the $\frac{1}{^{\alpha }}$-convexification of
$X=L^{p,r}(\log L)^{\beta }$ is a Banach space. On the other hand,
from (\ref{funda}) and since $\left( X^{(1/\alpha )}\right)
^{(\alpha )}=$ $X$, we obtain that given $0<s<1\ $and $0<q\leq
\infty $ the function $m_{X^{(1/\alpha )},s,\alpha ,q}$ defined in
Theorem \ref{embteo} verifies that
\begin{equation}
m_{X^{(1/\alpha )},s,\alpha ,q}(t)\simeq \left\{
\begin{array}{ll}
\int_{t}^{1}\left( \frac{z^{\frac{s}{Q}-\frac{1}{p}}}{(1+\ln \frac{1}{z}%
)^{\beta }}\right) ^{\frac{\alpha q}{q-\alpha }}\frac{dz}{z} & \text{if \ \ }%
\alpha <q<\infty , \\
\int_{t}^{1}\left( \frac{z^{\frac{s}{Q}-\frac{1}{p}}}{(1+\ln \frac{1}{z}%
)^{\beta }}\right) ^{\alpha }\frac{dz}{z} & \text{if \ \ }q=\infty , \\
\sup_{z\in \lbrack t,1)}\frac{z^{\frac{s}{Q}-\frac{1}{p}}}{(1+\ln \frac{1}{z}%
)^{\beta }} & \text{if \ \ }0<q\leq \alpha .
\end{array}
\right.  \label{mpalpha}
\end{equation}

Henceforth we shall assume that $\left( \Omega ,d,\mu \right) $ is
doubling with doubling constant $C_{\mu }$ and upper dimension $Q$
which satisfies the non-collapsing condition.

\begin{theorem}
\label{Teolorentzlog}Let $X=L^{p,r}(\log L)^{\beta }$ be a
Lorentz-Zygmund
space on $\Omega $ ($0<p<\infty ,$ $0<r\leq \infty $, $\beta \in \mathbb{R}$%
). Let $0<s<1,$ $0<q\leq \infty $. The embedding
\begin{equation}
\left\| f\right\| _{L^{\infty }}\preceq \left\| f\right\| _{\dot{B}%
_{L^{p,r}(\log L)^{\beta },q}^{s}}+\left\| f\right\| _{L^{\min
(1,p,r)}+L^{\infty }},\text{\ }  \label{infinfy}
\end{equation}
holds in the following situations:

\begin{enumerate}
\item  $s>\frac{Q}{p}$.

\item  $s=\frac{Q}{p},$ and
\begin{equation}
\left\{
\begin{array}{ll}
\beta >\frac{1}{\min (1,p,r)}-\frac{1}{q} & \text{if }\min
(1,p,r)<q\leq
\infty , \\
\beta \geq 0 & \text{if }0<q\leq \min (1,p,r).
\end{array}
\right.  \label{betas}
\end{equation}
\end{enumerate}
\end{theorem}

\begin{proof}
Given $X=L^{p,r}(\log L)^{\beta }$ pick $0<\alpha \leq 1$ satisfying
(\ref {elealpha}) and let $m_{X^{(1/\alpha )},\alpha ,q}$ the
function defined as in (\ref{mpalpha}). We will see that
$m_{X^{(1/\alpha )},s,\alpha ,q}(0)<\infty ,$ which by Theorem
\ref{embteo} implies that (\ref{infinfy}) holds.

\begin{enumerate}
\item  Case $s>\frac{Q}{p}.$ An elementary computation shows that $%
m_{X^{(1/\alpha )},s,\alpha ,q}(0)<0.$

\item  Case $s=\frac{Q}{p}.$ Notice that condition (\ref{betas}) implies
that $\beta \geq 0,$ then

\begin{enumerate}
\item  If $\min (1,p,r)=\min (1,r),$ then by (\ref{elealpha}) $\alpha =\min
(1,r)$. Besides, condition $\beta >\frac{1}{\min (1,r)}-\frac{1}{q}\ $%
implies $\frac{\beta \min (1,r)q}{q-\min (1,r)}>1$ if $\min
(1,r)<q<\infty $ and $\beta \min (1,r)>1$ if $q=\infty ,$ hence
$m_{X^{(1/\alpha )},s,\alpha ,q}(0)<\infty .$

Finally, if $0<q\leq \min (1,r)$ and $\beta \geq 0,$ then
\begin{equation*}
m_{X^{(1/\alpha )},s,\alpha ,q}(0)=\sup_{z\in \lbrack
0,1)}\frac{1}{(1+\ln \frac{1}{z})^{\beta }}<\infty .
\end{equation*}

\item  If $\min (1,p,r)=p,$ then given $\beta >\frac{1}{p}-\frac{1}{q},$ we
select $\alpha <p$ such that $\beta >\frac{1}{\alpha }-\frac{1}{q}>\frac{1}{p%
}-\frac{1}{q}.$ This choice of $\alpha $ implies that $\frac{\beta \alpha q}{%
q-\alpha }>1$ if $\alpha <q<\infty $ and $\beta \alpha >1$ if
$q=\infty ,$ therefore $m_{X^{(1/\alpha )},\alpha ,q}(0)<\infty .\ $

In case that $p\leq q\leq \infty $, from $\alpha <p,$ if follows
that
\begin{equation*}
m_{X^{(1/\alpha )},s,\alpha ,q}(0)=\sup_{z\in \lbrack
0,1)}\frac{1}{(1+\ln \frac{1}{z})^{\beta }}<\infty
\end{equation*}
when $0<q\leq \alpha <p$.
\end{enumerate}
\end{enumerate}
\end{proof}

\begin{theorem}
Let $X=L^{p,r}(\log L)^{\beta }$ be a Lorentz-Zygmund space on $\Omega $ ($%
0<p<\infty ,$ $0<r\leq \infty $, $\beta \in \mathbb{R}$). Let $0<s<1,$ $%
0<q\leq \infty $. The following embedding holds.

\begin{enumerate}
\item  If $s<\frac{Q}{p},$ then
\begin{equation*}
\left( \int_{0}^{1}\left[ t^{\frac{1}{p}-\frac{s}{Q}}\left( 1+\ln \frac{1}{t}%
\right) ^{\beta }f^{\ast }(t)\right] ^{q}\frac{dt}{t}\right)
^{1/q}\preceq \left\| f\right\| _{\dot{B}_{L^{p,r}(\log L)^{\beta
},q}^{s}}+\left\| f\right\| _{L^{\min (1,r,p)}+L^{\infty }}
\end{equation*}
\ and
\begin{equation*}
\sup_{0<t<1}t^{\frac{1}{p}-\frac{s}{Q}}(1+\ln \frac{1}{t})^{\beta
}f^{\ast }(t)\preceq \left\| f\right\| _{\dot{B}_{L^{p,r}(\log
L)^{\beta },\infty }^{s}}+\left\| f\right\| _{L^{\min
(1,r,p)}+L^{\infty }}.
\end{equation*}

\item  If $s=\frac{Q}{p}$ and
\begin{equation*}
\left\{
\begin{array}{ll}
\beta \leq \frac{1}{\min (1,p,r)}-\frac{1}{q} & \text{if }\min
(1,p,r)<q\leq
\infty , \\
\beta <0 & \text{if }0<q\leq \min (1,p,r),
\end{array}
\right.
\end{equation*}
then

\begin{enumerate}
\item  If $\min (1,r)<p,$ or $\min (1,r)=p$ an $\beta \geq 0,$ then

\begin{enumerate}
\item  If $\beta =\frac{1}{\min (1,r)}-\frac{1}{q}$ and $\min (1,r)<q\leq
\infty ,$ then
\begin{equation*}
\left( \int_{0}^{1}\frac{f^{\ast }(t)^{q}}{(1+\ln \frac{1}{t})\left(
1+\ln
\left( 1+\ln \frac{1}{t}\right) \right) ^{\frac{q}{\min (1,r)}}}\frac{dt}{t}%
\right) ^{1/q}\preceq \left\| f\right\| _{\dot{B}_{L^{p,r}(\log
L)^{\beta },q}^{s}}+\left\| f\right\| _{L^{\min (1,r)}+L^{\infty }},
\end{equation*}
if $\min (1,r)<q<\infty ,$ and
\begin{equation*}
\sup\limits_{t\in \lbrack 0,1]}\frac{f^{\ast }(t)}{\left( 1+\ln
\left( 1+\ln
\frac{1}{t}\right) \right) ^{\beta }}\preceq \left\| f\right\| _{\dot{B}%
_{L^{p,r}(\log L)^{\beta },\infty }^{s}}+\left\| f\right\| _{L^{\min
(1,r)}+L^{\infty }}.
\end{equation*}

\item  If $\beta <\frac{1}{\min (1,r)}-\frac{1}{q}$ and $\min (1,r)<q\leq
\infty ,$ then
\begin{equation*}
\left( \int_{0}^{1}\left( \frac{f^{\ast }(t)}{\left( 1+\ln \frac{1}{t}%
\right) ^{\frac{1}{\min (1,r)}-\beta }}\right)
^{q}\frac{dt}{t}\right) ^{1/q}\preceq \left\| f\right\|
_{\dot{B}_{L^{p,r}(\log L)^{\beta },q}^{s}}+\left\| f\right\|
_{L^{\min (1,r)}+L^{\infty }},
\end{equation*}
if $\min (1,r)<q<\infty ,$ and
\begin{equation*}
\sup\limits_{t\in \lbrack 0,1]}\frac{f^{\ast }(t)}{\left( 1+\ln \frac{1}{t}%
\right) ^{\frac{1}{\min (1,r)}-\beta }}\preceq \left\| f\right\| _{\dot{B}%
_{L^{p,r}(\log L)^{\beta },\infty }^{s}}+\left\| f\right\| _{L^{\min
(1,r)}+L^{\infty }}.
\end{equation*}

\item  If $\beta <0\ $and $0<q\leq \min (1,r),$ then
\begin{equation*}
\left( \int_{0}^{1}f^{\ast }(t)^{q}\left( 1+\ln \frac{1}{z}\right)
^{\beta
q-1}\frac{dt}{t}\right) ^{1/q}\preceq \left\| f\right\| _{\dot{B}%
_{L^{p,r}(\log L)^{\beta },q}^{s}}+\left\| f\right\| _{L^{\min
(1,r)}+L^{\infty }}.
\end{equation*}
\end{enumerate}

\item  If $p<\min (1,r)$ or $\min (1,r)=p$ and $\beta <0,$ then for any $%
\alpha <p$ we have that

\begin{enumerate}
\item  If $\beta \leq \frac{1}{p}-\frac{1}{q}\ $and $p\leq q\leq \infty ,$
then there is positive constant $C_{\alpha }$ that blows up when $a$
tends to $p,$ such that
\begin{equation*}
\left( \int_{0}^{1}\left( \frac{f^{\ast }(t)}{\left( 1+\ln \frac{1}{t}%
\right) ^{\frac{1}{\alpha }-\beta }}\right) ^{q}\frac{dt}{t}\right)
^{1/q}\leq C_{\alpha }\left( \left\| f\right\|
_{\dot{B}_{L^{p,r}(\log L)^{\beta },q}^{s}}+\left\| f\right\|
_{L^{p}+L^{\infty }}\right) ,
\end{equation*}
if $p\leq q<\infty ,$ and
\begin{equation*}
\sup_{t\in \lbrack 0,1]}\frac{f^{\ast }(t)}{(1+\ln \frac{1}{t})^{\beta -%
\frac{1}{\alpha }}}\leq C_{\alpha }\left( \left\| f\right\| _{\dot{B}%
_{L^{p,r}(\log L)^{\beta },\infty }^{s}}+\left\| f\right\|
_{L^{p}+L^{\infty }}\right)
\end{equation*}

\item  If $\beta <0$ and $0<q<p,$ then
\begin{equation*}
\left( \int_{0}^{1}f^{\ast }(t)^{q}\left( 1+\ln \frac{1}{z}\right)
^{\beta
q-1}\frac{dt}{t}\right) ^{1/q}\preceq \left\| f\right\| _{\dot{B}%
_{L^{p,r}(\log L)^{\beta },q}^{s}}+\left\| f\right\|
_{L^{p}+L^{\infty }}.
\end{equation*}
\end{enumerate}
\end{enumerate}
\end{enumerate}
\end{theorem}

\begin{proof}
Let $0<\alpha \leq 1$ satisfying (\ref{elealpha}), $0<s<1$ and
$0<q\leq \infty .$ Let $m(t)=m_{X^{(1/\alpha )},s,\alpha ,q}(t)$ the
function defined as in (\ref{mpalpha}). It is a matter of a tedious
but elementary calculation to verify that if $s<\frac{Q}{p},$ then
$m(0)=\infty $ and
\begin{equation*}
\left\{
\begin{array}{ll}
\frac{\left| m^{\prime }(t)\right| }{\left( 1+m(t)\right) ^{q/\alpha }}%
\simeq \left( t^{\frac{1}{p}-\frac{s}{Q}}(1+\ln \frac{1}{t})^{\beta
}\right)
^{q}\frac{1}{t} & \text{if \ \ }\alpha <q<\infty , \\
\frac{1}{\left( 1+m(t)\right) ^{1/\alpha }}\simeq t^{\frac{1}{p}-\frac{s}{Q}%
}(1+\ln \frac{1}{t})^{\beta } & \text{if \ \ }q=\infty , \\
\int_{0}^{t}v^{q}(z)dz\preceq v^{q}(t) & \text{if \ }0<q\leq \alpha
,
\end{array}
\right.
\end{equation*}
consequently part 1 of Theorem \ref{embteo} applies.

In the same way, if $s=\frac{Q}{p},$ then $m(0)=\infty ,$ and

\begin{equation*}
\left\{
\begin{array}{ll}
\frac{\left| m^{\prime }(t)\right| }{\left( 1+m(t)\right) ^{q/\alpha }}%
\simeq \frac{1}{t(1+\ln \frac{1}{t})}\frac{1}{\left( 1+\ln \left(
1+\ln \frac{1}{t}\right) \right) ^{q/\alpha }}, & \text{if }\alpha
<q<\infty \text{
and }\frac{\beta \alpha q}{q-\alpha }=1, \\
\frac{\left| m^{\prime }(t)\right| }{\left( 1+m(t)\right) ^{q/\alpha }}%
\simeq \frac{1}{t}(1+\ln \frac{1}{t})^{\left( \beta \alpha -1\right) \frac{q%
}{\alpha }}, & \text{if }\alpha <q<\infty \text{ and }\frac{\beta \alpha q}{%
q-\alpha }<1, \\
\frac{1}{\left( 1+m(t)\right) ^{1/\alpha }}\simeq \frac{1}{\left(
1+\ln \left( 1+\ln \frac{1}{t}\right) \right) ^{\beta }}, & \text{if
}q=\infty
\text{ and }\alpha \beta =1, \\
\frac{1}{\left( 1+m(t)\right) ^{1/\alpha }}\simeq (1+\ln
\frac{1}{t})^{\beta
-\frac{1}{\alpha }}, & \text{if }q=\infty \text{ and }\alpha \beta <1, \\
\int_{0}^{t}\left( \left( 1+\ln \frac{1}{z}\right) ^{\beta -\frac{1}{q}%
}\right) ^{q}\frac{dz}{z}\preceq \left( 1+\ln \frac{1}{t}\right)
^{\beta q}, & \text{if }0<q\leq \alpha .
\end{array}
\right.
\end{equation*}
Thus:

\begin{enumerate}
\item  If $\min (1,r)<p,$ or $\min (1,r)=p$ and $\beta \geq 0,$ then\ $%
\alpha =\min (1,r)$ and the part $(a)$ of statement $2$ follows by
Theorem \ref{embteo}.

\item  If $p<\min (1,r)$ or $\min (1,r)=p$ and $\beta <0,$ then for any $%
\alpha <p,$ we have that $\beta \leq
\frac{1}{p}-\frac{1}{q}<\frac{1}{\alpha }-\frac{1}{q}$, therefore
using one more time Theorem \ref{embteo} we obtain the $(b)$-part of
statement $2$.
\end{enumerate}
\end{proof}

As a byproduct of the above Theorem, having $\beta =0,$ we obtain
the
following embedding result for Lorentz-Haj\l asz-Besov spaces $\dot{B}%
_{L^{p,r},q}^{s}$.

\begin{corollary}
Let $0<p<\infty ,$ $0<r\leq \infty ,$ $0<q\leq \infty $ and $0<s<1.$
The following assertions are true:

\begin{enumerate}
\item  If $s<\frac{Q}{p}$, then
\begin{equation*}
\left( \int_{0}^{1}\left[ t^{\frac{1}{p}-\frac{s}{Q}}f^{\ast }(t)\right] ^{q}%
\frac{dt}{t}\right) ^{1/q}\preceq \left\| f\right\| _{\dot{B}%
_{L^{p,r},q}^{s}}+\left\| f\right\| _{L^{\min (1,p,r)}+L^{\infty }},
\end{equation*}
if $0<q<\infty $, and
\begin{equation*}
\sup_{0<t<1}t^{\frac{1}{p}-\frac{s}{Q}}f^{\ast }(t)\preceq \left\|
f\right\| _{\dot{B}_{L^{p,r},\infty }^{s}}+\left\| f\right\|
_{L^{\min (1,p,r)}+L^{\infty }}.
\end{equation*}

\item  If $s=\frac{Q}{p}$, then

\begin{enumerate}
\item  if $0<r<p,$ then
\begin{equation*}
\left( \int_{0}^{1}\left[ \frac{f^{\ast }(t)}{\left( 1+\ln \frac{1}{t}%
\right) ^{\frac{1}{\min (1,r)}}}\right] ^{q}\frac{dt}{t}\right)
^{1/q}\preceq \left\| f\right\| _{\dot{B}_{L^{p,r},q}^{s}}+\left\|
f\right\| _{L^{\min (1,r)}+L^{\infty }}\text{,}
\end{equation*}
if $\min (1,r)<q<\infty $, and
\begin{equation*}
\sup_{0<t<1}\frac{f^{\ast }(t)}{\left( 1+\ln \frac{1}{t}\right) ^{\frac{1}{%
\min (1,r)}}}\preceq \left\| f\right\| _{\dot{B}_{L^{p,r},\infty
}^{s}}+\left\| f\right\| _{L^{\min (1,p)}+L^{\infty }}.
\end{equation*}

\item  if $1<p\leq r,$ or $1\leq p=r,$ then
\begin{equation*}
\left( \int_{0}^{1}\left[ \frac{f^{\ast }(t)}{\left( 1+\ln \frac{1}{t}%
\right) }\right] ^{q}\frac{dt}{t}\right) ^{1/q}\preceq \left\| f\right\| _{%
\dot{B}_{L^{p,r},q}^{s}}+\left\| f\right\| _{L^{1}+L^{\infty
}}\text{,}
\end{equation*}
if $1<q<\infty $, and
\begin{equation*}
\sup_{0<t<1}\frac{f^{\ast }(t)}{1+\ln \frac{1}{t}}\preceq \left\|
f\right\| _{\dot{B}_{L^{p,r},\infty }^{s}}+\left\| f\right\|
_{L^{\min (1,p)}+L^{\infty }}.
\end{equation*}

\item  if $0<p\leq r$ and $0<p\leq 1,$ then for any $\alpha <p$
\begin{equation*}
\left( \int_{0}^{1}\left[ \frac{f^{\ast }(t)}{\left( 1+\ln \frac{1}{t}%
\right) ^{\frac{1}{\alpha }}}\right] ^{q}\frac{dt}{t}\right)
^{1/q}\preceq \left\| f\right\| _{\dot{B}_{L^{p,r},q}^{s}}+\left\|
f\right\| _{L^{\alpha }+L^{\infty }}\text{,}
\end{equation*}
and
\begin{equation*}
\sup_{0<t<1}\frac{f^{\ast }(t)}{\left( 1+\ln \frac{1}{t}\right) ^{\frac{1}{%
\alpha }}}\preceq \left\| f\right\| _{\dot{B}_{L^{p,r},\infty
}^{s}}+\left\| f\right\| _{L^{\alpha }+L^{\infty }}.
\end{equation*}
\end{enumerate}

\item  If $s>\frac{Q}{p}$ or $s=\frac{Q}{p}$ and $q<p,$ then
\begin{equation*}
\left\| f\right\| _{L^{\infty }}\preceq \left\| f\right\| _{_{\dot{B}%
_{L^{p,r},q}^{s}}}+\left\| f\right\| _{L^{\alpha }+L^{\infty }}.
\end{equation*}
\end{enumerate}
\end{corollary}

Similarly, if $p=r$ then we obtain embeddings for
Lorentz-Zygmund-Haj\l asz-Besov spaces $B_{L^{p}(\log L)^{\beta
}}^{s}$, we left the details to the interested reader.


\begin{thebibliography}{99}
\bibitem{AGH}  R. Alvarado, P. G\'{o}rka, and P. Haj{\l }asz, \textsl{%
Sobolev embedding for $M^{1,p}$ spaces is equivalent to a lower
bound of the measure}. J. Funct. Anal. \textbf{279} (2020), 1--39.

\bibitem{AYY}
R. Alvarado, D. Yang and W. Yuan,  \textsl{Optimal embeddings for Triebel–Lizorkin and Besov spaces on quasi-metric measure spaces},  arXiv: 2202.06389

\bibitem{AN} S. V. Astashkin and P. G. Nilsson,  \textsl{A description of interpolation spaces for quasi-Banach couples by real
K-method}, arXiv:2112.13248.


\bibitem{BL}  J. Bergh and J. L\"{o}fstr\"{o}m, \textsl{Interpolation
spaces. An introduction}. Grundlehren der Mathematischen
Wissenschaften, No. 223. SpringerVerlag, Berlin New York, 1976.

\bibitem{BHS}  J. Bastero, H. Hudzik, and A.M. Steinberg, \textsl{On
smallest and largest spaces among rearrangement-invariant $p$-Banach
function spaces $(0<p<1)$} Indag. Math. N.S. \textbf{2}(3) (1991),
283--288.

\bibitem{BeRu}  C. Bennett and K. Rudnick, On LorentzZygmund spaces,
Dissert. Math. \textbf{175} (1980), 1--72.

\bibitem{BS}  C. Bennett and R. Sharpley, \textsl{Interpolation of Operators}%
, Academic Press, Boston\textbf{\ }(1988).

\bibitem{BreWa}  H. Brezis and S. Wainger, \textsl{A note on limiting cases
of Sobolev embeddings and convolution inequalities,} Comm. Partial
Di . Eq. \textbf{5} (1980), 773--789.

\bibitem{CaoGri}  J.Cao and A. Grigory'an, \textsl{Heat Kernels and Besov
Spaces on Metric Measure Spaces.} Preprint

\bibitem{CosMi}  S. Costea and M. Miranda Jr, \textsl{Newtonian Lorentz
metric spaces}, Illinois J. Math. \textbf{56} (2012), no. 2,
579--616.

\bibitem{Curbera}  G.P. Curbera, J. Garc\'{i}a-Cuerva, J.M. Martell, and C.
P\'{e}rez, \textsl{Extrapolation with weights,
rearrangement-invariant
function spaces, modular inequalities and applications to singular integrals.%
} Adv. Math. \textbf{203} (2006), no. 1, 256--318.



\bibitem{Amiran}  A. Gogatishvili, P. Koskela, and N. Shanmugalingam,
\textsl{Interpolation properties of Besov spaces defined on metric
spaces}, Math. Nachr. \textbf{283} (2010), 215--231.

\bibitem{Amiran1}  A. Gogatishvili, P. Koskela, and Y. Zhou, \textsl{%
Characterizations of Besov and Triebel-Lizorkin Spaces on Metric
Measure Spaces,} Forum Math. \textbf{25} (2013), 787--819.

\bibitem{Gor}  P. G\'{o}rka, \textsl{In metric-measure spaces Sobolev
embedding is equivalent to a lower bound for the measure.} Potential
Anal. \textbf{47} (2017), no. 1, 13--19.

\bibitem{GK}  L. Grafakos and N.J. Kalton, Some remarks on multilinear maps
and interpolation, Math. Ann. \textbf{319} (2001), 151--180.

\bibitem{HMY}  Y. Han, D. M\"{u}ller, and D. Yang, \textsl{A theory of Besov
and Triebel-Lizorkin spaces on metric measure spaces modeled on
Carnot-Caratheodory spaces}, Abstr. Appl. Anal. \textbf{2008}, Art.
ID 893409, 250 pp.

\bibitem{Ha2}  P. Haj{\l }asz, \textsl{\ Sobolev spaces on an arbitrary
metric spaces}, Potential Anal. \textbf{5} (1996), 403--415.

\bibitem{Ha1}  P. Haj{\l }asz, \textsl{Sobolev spaces on metric-measure
spaces. In Heat kernels and analysis on manifolds, graphs, and
metric spaces (Paris, 2002)}, Contemp. Math, Amer. Math. Soc.
Providence, RI. \textbf{338} (2003), 173--218.

\bibitem{HIH}  T. Heikkinen, L. Ihnatsyeva, and H. Tuominen, \textsl{Measure
density and extension of Besov and Triebel-Lizorkin functions,} J.
Fourier Anal. Appl. \textbf{22} (2016), no. 2, 334--382.

\bibitem{HKa}T. Heikkinen and  N. Karak,  \textsl{Orlicz-Sobolev embeddings, extensions and Orlicz-PoincarÃ© inequalities,} J. Funct. Anal. \textbf{282}(2) (2022), 109292.

\bibitem{JS}  W.B. Johnson and G. Schechtman, \textsl{Sums of independent
random variables in rearrangement invariant function spaces,} Ann.
Probab. \textbf{17} (1989) 789--808

\bibitem{Kal}  N.J. Kalton,\textsl{\ Convexity conditions on non-locally
convex lattices,} Glasgow J. Math. \textbf{25} (1984), 141--152.

\bibitem{Kal1}  N.J. Kalton,\textsl{Plurisubharmonic functions on quasi-Banach spaces}, Studia Math.
\textbf{84} (1987) 297--324


\bibitem{Ka01}  N. Karak,\textsl{\ Lower bound of measure and embeddings of
Sobolev, Besov and Triebel-Lizorkin spaces,} Math. Nachr.
\textbf{293} (2020), no. 1, 120--128.

\bibitem{ko}  V. I. Kolyada, \textsl{Estimates of rearrangements and
embedding theorems}, Mat. Sb. \textbf{136 } (1988), 3--23 (in
Russian), English transl.: Math. USSR-Sb. \textbf{55 }(1989), 1--21
.

\bibitem{KoSak}  P. Koskela and E. Saksman, \textsl{Pointwise
characterizations of Hardy-Sobolev functions}, Math. Res. Lett.
\textbf{15} (2008), no. 4, 727--744.

\bibitem{KPS}  S.G. Krein, Yu.I. Petunin, and E.M. Semenov, \textsl{%
Interpolation of Linear Operators, } Transl. Math. Monogr. Amer.
Math. Soc. \textbf{54} (1982).

\bibitem{LT}  J. Lindenstrauss and L. Tzafriri, Classical Banach spaces. II,
Springer-Verlag, Berlin-New York, 1979.

\bibitem{Mal2}  L. Mal\'{y}, \textsl{Minimal weak upper gradients in
Newtonian spaces based on quasi-Banach function lattices,} Ann.
Acad. Sci. Fenn. Math. \textbf{38} (2013), no. 2, 727--745.

\bibitem{Mal1}  L. Mal\'{y}, \textsl{Newtonian spaces based on quasi-Banach
function lattices,} Math. Scand. \textbf{119} (2016), no. 1,
133--160.

\bibitem{Mar}  J. Mart\'{i}n, \textsl{Symmetrization inequalities in the
fractional case and Besov embeddings}, J. Math. Anal. Appl.
\textbf{344} (2008), 99--123.

\bibitem{MMilPAMS}  J. Mart\'{i}n and M. Milman, \textsl{Symmetrization
inequalities and Sobolev embeddings}, Proc. Amer. Math. Soc.
\textbf{134} (2006), no. 8, 2335--2347.

\bibitem{MM6}  J. Mart\'{i}n and M. Milman, \textsl{Isoperimetry and
Symmetrization for Logarithmic Sobolev inequalities}\emph{, }J.
Funct. Anal. \textbf{256} (2009), 149--178.

\bibitem{MM3}  J. Mart\'{\i}n and M. Milman, \textsl{Pointwise
symmetrization inequalities for Sobolev functions and applications},
Adv. Math. \textbf{225} (2010), 121--199.

\bibitem{MM4}  J. Mart\'{\i}n and M. Milman, \textsl{Fractional Sobolev
inequalities: symmetrization, isoperimetry and interpolation,}
Ast\'{e}risque. \textbf{366} (2014), x+127 pp.

\bibitem{MaOr}  J. Mart\'{i}n and W.A. Ortiz, \textsl{A Sobolev type
embedding theorem for Besov spaces defined on doubling metric
spaces.} J. Math. Anal. Appl. \textbf{479} (2019), no. 2,
2302--2337.

\bibitem{MaOr1}  J. Mart\'in, and W.A. Ortiz, \textsl{Sobolev embeddings for
Fractional Haj{\l}asz-Sobolev spaces in the setting of rearrangement
invariant spaces}, Potential Analysis (to appear).

\bibitem{Mas}  M. Masty{\l}o, \textsl{The modulus of smoothness in metric
spaces and related problems}, Potential Anal. \textbf{35} (2011)
301--328.

\bibitem{MO} W. Matuszewska and W. Orlicz, \textsl{A note on the theory of s-normed spaces
of f-integrable functions}, Studia Math. \textbf{21} (1961/62) 107--115.

\bibitem{MY}  D. M\"{u}ller, and D. Yang, \textsl{A difference
characterization of Besov and Triebel-Lizorkin spaces on RD-spaces},
Forum Math. \textbf{21} (2009), 259--298.

\bibitem{Mu}  T. Muramatu, \textsl{On Besov spaces and Sobolev spaces of
generalized functions defined on a general region}, Publ. Res. Inst.
Math. Sci. \textbf{9} (1973/74), 325--396.

\bibitem{Sawyer}  E. Sawyer, \textsl{Boundedness of classical operators on
classical Lorentz spaces}, Studia Math. \textbf{96} (1990),
145--158.

\bibitem{Tuo}  H. Tuominen, \textsl{Orlicz-Sobolev spaces on metric measure
spaces}, Ann. Acad. Sci. Fenn. Math. Dissertationes 135.

\bibitem{YZ}  D. Yang and Y. Zhou, \textsl{New properties of Besov and
Triebel-Lizorkin spaces on RD-spaces}, Manuscripta Math.
\textbf{134} (2011), 59--90.
\end{thebibliography}
\end{document}